\documentclass{siamart250106}

\usepackage{graphicx,enumitem}
\usepackage{amsfonts, amssymb}
\usepackage{amsmath,mathrsfs}
\usepackage{cancel,overpic,MnSymbol}
\usepackage{afterpage,xcolor}
\numberwithin{equation}{section}

\newcommand{\R}{\mathbb{R}}
\newcommand{\ddx}[2]{\partial_{#2}#1}
\newcommand{\dfdx}[2]{\frac{\partial #1}{\partial #2}}
\renewcommand{\phi}{\varphi}

\newcommand{\F}{\mathcal{F}}

\newcommand{\rx}{\tilde{r}^\times}%

\newcommand{\B}{\mathcal{B}}

\newcommand{\tir}{\tilde{r}} 
\newcommand{\tix}{\tilde{x}}
\newcommand{\tip}{\tilde{p}}
\newcommand{\tik}{\nu}       

\newcommand{\LL}{\mathscr{L}}

\newtheorem{thrm}[theorem]{Theorem}
\newtheorem{lmm}[theorem]{Lemma}
\theoremstyle{definition}

\newtheorem{rmrk}[theorem]{Remark}

\newcommand{\tipm}{\tilde{P}}

\title{A First Closed-Form Robin Parameter for Optimized Schwarz Waveform
Relaxation as a Function of Iteration Count, Overlap, and Time Window Length}

\author{Alejandro {Alfonso Rodriguez}\footnotemark[1] \and Martin J.~Gander\footnotemark[2] \and Felix Kwok\footnotemark[3]}

\begin{document}
\maketitle
\renewcommand{\thefootnote}{\fnsymbol{footnote}}

\begin{abstract}Schwarz waveform relaxation (SWR) is a parallel-in-time iterative
method for the solution of evolution problems and belongs to the class
of Parallel-in-Time algorithms that have received increasing
attention over the past three decades. For parabolic problems, SWR
converges linearly over long time intervals and exhibits superlinear
convergence on short time intervals. Since the introduction of
optimized variants (OSWR) nearly three decades ago, determining the
optimized Robin parameter when OSWR is applied on bounded time windows
has remained an open problem.

We introduce here a new analytical approach that, for the first time,
yields a closed-form expression for the optimal Robin parameter in the
transmission conditions for the one-dimensional heat equation with two
overlapping subdomains. The resulting formula explicitly characterizes
the dependence of the optimal parameter on the iteration number, the
overlap size, and the length of the time window. Beyond resolving this
long-standing optimization problem, our analysis provides a framework
that can be extended to more general settings, including
higher-dimensional problems, multiple subdomains, and higher-order
transmission conditions.
\end{abstract}

\pagestyle{myheadings}
\thispagestyle{plain}
\markboth{A.~ALFONSO RODRIGUEZ, M.~J. GANDER AND F. KWOK}{A closed-form Robin parameter for OSWR}

\footnotetext[1]{Department of Mathematics and Statistics, Universit\'e Laval, Quebec, Canada {\tt alejandro.alfonso-rodriguez.1@ulaval.ca}.}
\footnotetext[2]{Section de math\'ematiques, Universit\'e de Gen\`eve, Geneva, Switzerland {\tt martin.gander@unige.ch}.}
\footnotetext[3]{Department of Mathematics and Statistics, Universit\'e Laval, Quebec, Canada {\tt felix.kwok@mat.ulaval.ca}.}

\renewcommand{\thefootnote}{\arabic{footnote}}

\section{Introduction}
The convergence behaviour of classical and optimized Schwarz Waveform Relaxation (SWR) methods for parabolic partial differential equations (PDEs) 
differs markedly from that of Schwarz methods for elliptic problems \cite{GanderOSM}, owing to the non-normal nature of the time derivative operator. 
It is known that classical SWR methods (i.e., with Dirichlet transmission conditions) converge superlinearly on bounded time intervals, but such convergence deteriorates as the time interval
length $T$ goes to infinity; in the limit, convergence is only linear \cite{GanderStuart,GiladiKeller,GanderZhao}. 
As the iteration progresses, the convergence curve exhibits different phases: there is an initial rapid phase, then a slowdown, and finally an acceleration to superlinear convergence
\cite{GanderMartin}. Moreover, when Robin transmission conditions are used,
the best Robin parameter $p$ that minimizes error reduction depends on the number of iterations $n$ one intends to take, rather than a fixed constant \cite{Courvoisier}.

The classical approach to finding best parameters via Fourier analysis, see \cite{Bennequin}, is unable to capture this dependence on $n$. Indeed, such an analysis does
not even yield superlinear convergence estimates (but see \cite{GanderMartin} for a Kernel Fourier analysis that uses the heat kernel to restore superlinearity in the estimates).
As a result, superlinear estimates typically employ kernel techniques \cite{GanderStuart, GanderZhao} or Laplace transforms \cite{Mandal, Mandal2, Mandal3}. The latter is
useful when the iteration function in Laplace space admits a closed form inverse Laplace transform, or can be expressed as a convergent series of such inverse transforms. 
This is feasible for Dirichlet and Neumann transmission conditions, but not for Robin transmission conditions. As a result, the problem of finding $n$-dependent best Robin parameters
remains so far an open problem, even for the simple case of the one-dimensional heat equation in two subdomains. Note that in the discrete-time setting, it is possible to find the
best Robin parameter by numerically minimizing the norm of the appropriate iteration matrix, whose size is proportional to the number of time steps, see \cite{arnoult2023discrete,arnoult2025methode}.

The goal of this paper is twofold. The first goal is to introduce a theoretical tool for estimating the error of domain decomposition methods with simple Laplace transforms,
even when the inverse transforms do not have closed form formulas. The tool relies on a simple relationship between the Laplace and Fourier transforms, together with the 
Plancherel theorem. Although the tool is easy to understand, its application is not so widespread in the domain decomposition literature; we therefore illustrate its utility by
giving a new convergence analysis of the classical SWR method for two subdomains, for which we are able to recover the known superlinear estimates using only a few lines of calculations. The second goal of the paper is to use the same techniques to fully analyze the convergence of optimized SWR (OSWR), first as a function of the Robin parameter $p$, and then
by choosing the parameter $p$ to minimize this contraction estimate for a given number of iterations $n$. 
Note that this is the first time we are able to truly optimize $p$ in the superlinear regime, even though the algorithm exists and has been used for three decades.
In each case, we use numerical examples to illustrate the sharpness of our estimates.

\section{Laplace transform and exponential weighting}
Let $\Omega = \R$. To explain this new technique, it suffices to focus on the one-dimensional heat equation\footnote{We use $y$ to denote the spatial variable, since $x$ will be used for another quantity later in the manuscript.}
\begin{equation}\label{eq:1dheat_with_f}
\ddx{u}{t} - \ddx{u}{yy} = f,\qquad (y,t)\in \Omega \times (0,T)
\end{equation}
with source term $f\in L^2(0,T; L^2(\Omega))$ and initial conditions $u(\cdot,0) = u_0 \in H^1(\Omega)$. 
Note the higher regularity assumed in the initial conditions, which is needed for the well-posedness of the SWR methods later. 
The problem \eqref{eq:1dheat_with_f} is well posed \cite{GanderHalpernCRAS,GanderHalpernSINUM2007} with a unique solution $u\in H^{2,1}(\Omega\times (0,T)):= L^2(0,T; H^2(\Omega)) \cap H^1(0,T; L^2(\Omega))$. Next,
we consider a decomposition of $\Omega$ into two spatial subdomains
$\Omega_1 = (-\infty, \delta)$, $\Omega_2 = (0,\infty)$ with an overlap of width $\delta > 0$. The alternating (O)SWR method solves, for $n=1,2,3,\ldots$,
\begin{alignat*}{2}
\ddx{u_1^n}{t} - \ddx{u_1^n}{yy}  &= f~~~~~~~~~~~~~~~~\text{in $\Omega_1\times (0,T)$,}\qquad &
 \ddx{u_2^n}{t} - \ddx{u_2^n}{yy} &= f~~~~~~~~~~~~~\text{in $\Omega_2\times (0,T)$,}\\
u_1^n(\cdot,0) &= u_0~~~~~~~~~~~~~~~\text{in $\Omega_1$}, & 
u_2^n(\cdot,0) &= u_0~~~~~~~~~~~~\text{in $\Omega_2$},\\
\B_1u_1^n(\delta,\cdot) &= \B_1u_2^{n-1}(\delta,\cdot)\quad\text{on $(0,T)$}, & 
\B_2u_2^n(0,\cdot) &= \B_2u_1^{n}(0,\cdot)\quad\text{on $(0,T)$,}
\end{alignat*}
where the transmission conditions can be
\begin{itemize}
\item classical (Dirichlet) conditions: $\B_i = \mathrm{id}$, i.e., the identity operator;
\item optimized (Robin) conditions: $\B_i = \partial_{n_i} + p$, where $p > 0$ is a constant and $\partial_{n_i}$ is the normal derivative, with $n_i$ being the normal vector pointing outward from $\Omega_i$. In other words, we have $\B_1 = \partial_y + p$ and $\B_2 = -\partial_y + p$.
\end{itemize}
For classical SWR, the iteration is well posed if $f\in L^2(0,T; L^2(\Omega))$, $u_0|_{\Omega_i} \in H^1(\Omega_i)$, and if the initial guess $u_1^1(\delta,\cdot)$
is in $H^{3/4}(0,T)$ and compatible with $u_0$ at $t=0$. The iterates $u^n_i$ then converge to $u|_{\Omega_i}$ in $H^{2,1}(\Omega_i\times (0,T))$ as $n\to \infty$. For
OSWR with $p>0$, the algorithm is well-posed if  $f\in L^2(0,T; L^2(\Omega))$, $u_0|_{\Omega_i} \in H^1(\Omega_i)$ and if the initial Robin trace $\B_1u_1^1(\delta,\cdot)$,
is in $H^{1/4}(0,T)$; in that case, the iterates $u_i^n$ all lie in $H^{2,1}(\Omega_i\times (0,T))$ and converge to $u|_{\Omega_i}$ as
$n\to \infty$. See \cite{GanderHalpernCRAS,GanderHalpernSINUM2007} for the proof of these results.

Since the problem \eqref{eq:1dheat_with_f} is linear, it suffices to analyze the convergence of the local errors $E_i^n = u_i^n - u|_{\Omega_i}$ to zero. The error equations
satisfied by $E_i^n$ have the same form as the original problem, but with zero initial conditions and source term:
\begin{equation}
\label{err_eq}
\begin{aligned}
\ddx{E_1^n}{t} - \ddx{^2E_1^n}{y}  &= 0~~~~~~~~~~~~~~~~\text{in $\Omega_1\times (0,T)$,}\qquad & 
\ddx{E_2^n}{t} - \ddx{^2E_2^n}{y} &= 0~~~~~~~~~~~~~\text{in $\Omega_2\times (0,T)$,}\\
E_1^n(\cdot,0) &= 0~~~~~~~~~~~~~~~~\text{in $\Omega_1$}, & 
E_2^n(\cdot,0) &= 0~~~~~~~~~~~~~\text{in $\Omega_2$},\\
\B_1E_1^n(\delta,\cdot) &= \B_1E_2^{n-1}(\delta,\cdot)\quad\text{on $(0,T)$}, & 
\B_2E_2^n(0,\cdot) &= \B_2E_1^{n}(0,\cdot)\quad\text{on $(0,T)$.}
\end{aligned}
\end{equation}

We now introduce the main idea of our paper, namely the Laplace transform as an exponentially weighted Fourier transform. Let 
$g \in L^2(0,T)$, and $\tilde{g} \in L^2(\R)$ be the extension of $g$ by zero outside the interval $(0,T)$, i.e.,
$$ \tilde{g}(t) = \begin{cases}
g(t), & 0 < t < T,\\
0, & \text{else}.
\end{cases}$$
Recall that the Laplace transform of $g$, which is a function of the complex variable $s = r+i\omega$ with $r \geq 0$, is given by
$$ \LL\{g\}(s) = \int_0^\infty g(t) e^{-st}\,dt = \int_{-\infty}^\infty e^{-i\omega t} (e^{-rt}\tilde{g}(t))\,dt = \F\{e^{-rt}\tilde{g}(t)\} =: \hat{g}(r,\omega), $$
where $\F$ denotes the Fourier transform. We can therefore invoke the Plancherel theorem to deduce that
\begin{equation}
\label{plancherel}
\frac{1}{2\pi}\int_{-\infty}^\infty |\hat{g}(r,\omega)|^2\,d\omega = \int_{-\infty}^\infty |e^{-rt}\tilde{g}(t)|^2\,dt = \int_0^T e^{-2rt}g^2(t)\,dt,
\end{equation}
which is valid for all $r\geq 0$.

We now explain how this idea can be used to analyze (O)SWR methods. To illustrate the technique, let us concentrate for now on the classical transmission
conditions $\B_1 = \B_2 = \mathrm{id}$; the case of optimized transmission conditions will be considered in detail in the next sections. We start by considering the
error equations \eqref{err_eq} and study the behaviour of $E_i^n(y,t)$ for all time $t>0$ (instead of truncating the solution at $t = T$). 
Writing $\hat{E}_1^n(y,s) = \LL\{E_1^n(y,\cdot)\}(s)$, the Laplace transform of the error equation for each $y\in \Omega_1$ reads
$$ s\hat{E}_1^n - \partial_{yy}E_{1}^n = 0, \qquad y \in (-\infty,\delta),$$
with $\lim_{y\to-\infty}\hat{E}_1(y,s) = 0$ and $\hat{E}^{n}_1(\delta,s) = \hat{E}_2^{n-1}(0,s)$ for all $\Re(s) \geq 0$. For a fixed $s$, this ODE in $y$ can be
solved explicitly to yield
\begin{equation}
\label{eq:recur_1}
 \hat{E}_1^n(y,s) = \hat{E}_2^{n-1}(\delta,s) e^{(y-\delta)\sqrt{s}},\qquad y < \delta, \quad n \geq 1,
\end{equation}
where $\sqrt{s} = \xi + i\eta$ is the unique square root of $s$ with positive real part. 
Similarly, taking the Laplace transform in time of the PDE in $\Omega_2$ and solving the resulting ODE gives
\begin{equation}
\label{eq:recur_2}
 \hat{E}_2^n(y,s) = \hat{E}_1^n(0,s)e^{-y\sqrt{s}}, \qquad y > 0, \quad n \geq 1.
\end{equation}
Combining \eqref{eq:recur_1} and \eqref{eq:recur_2} yields the recurrence
\begin{equation}
\label{eq:recur_3}
 \hat{E}_2^n(\delta,s) = \hat{E}_2^{n-1}(\delta,s)e^{-2\delta\sqrt{s}} = \cdots = \hat{E}_2^0(\delta,s)e^{-2n\delta\sqrt{s}}. 
\end{equation}
At this point, if we know the inverse Laplace transform $\rho_n(t)$ of $e^{-2n\delta\sqrt{s}}$, we can take inverse Laplace transforms of both sides to conclude that
$$ E_2^n(\delta,t) = \int_0^t \rho_n(t-\tau) E_2^0(\delta,\tau)\,d\tau, $$
and can use the properties of $\rho_n$ to estimate the size of $\|E_2^n(\delta,t)\|_{L^2(0,T)}$. Here, we take a different approach: we instead assume, without loss
of generality, that the initial error on the interface $E_2^0(\delta,\cdot)$ has support in the interval $(0,T)$.\footnote{By causality, any non-zero value of $E_2^0(\delta,t)$ at
$t > T$ will not impact the error of subsequent iterations inside the time interval $(0,T)$.}  Then for any $r \geq 0$, we write
\begin{align*}
 \|E_2^n(\delta,\cdot)\|_{L^2(0,T)}^2 = \int_0^T (E_2^n(\delta,t))^2\,dt &\leq e^{2rT}\int_{0}^{T} e^{-2rt} (E_2^n(\delta,t))^2\,dt \\
&\leq e^{2rT}\int_0^{{\infty}} e^{-2rt} (E_2^n(\delta,t))^2\,dt = \frac{e^{2rT}}{2\pi}\int_{-\infty}^\infty |\hat{E}_2^n(\delta,r + i\omega)|^2\,d\omega,
\end{align*}
where we used the Plancherel theorem \eqref{plancherel} twice.
With the recurrence \eqref{eq:recur_3}, we deduce that
\begin{align*}
\|E_2^n(\delta,\cdot)\|_{L^2(0,T)}^2 &\leq \frac{e^{2rT}}{2\pi}\int_{-\infty}^\infty |\hat{E}_2^0(\delta,r+i\omega)|^2|e^{-4n\delta\sqrt{r+i\omega}}|\,d\omega\\
&\leq \frac{(\rho(r,n))^2}{2\pi}\int_{-\infty}^\infty |\hat{E}_2^0(\delta,r+i\omega)|^2\,d\omega =(\rho(r,n))^2 \|E_2^0(\delta,\cdot)\|_{L^2(0,T)}^2,
\end{align*}
where we introduced
\begin{equation}
\label{eq:classical_factor}
\rho(r,n) = e^{rT} \sup_{\omega\in\R}|e^{-2n\delta\sqrt{r+i\omega}}|,
\end{equation}
an upper bound for the $n$-step contraction factor for the $L^2$-norm (without the square). Note that this upper bound is valid for any $r\geq 0$; thus, 
we can obtain the best estimate by choosing $r$ judiciously.
Indeed, if we let $\sqrt{r+i\omega} = x + ib$ and calculate
\begin{equation}\label{eq:b}
 r + i\omega = (x + ib)^2 = (x^2-b^2) + 2ibx \implies r = x^2-b^2,
\end{equation}
then together with $|r+i\omega| = |x + ib|^2 = x^2 + b^2$, we deduce that $x^2 = \frac12(r + \sqrt{r^2+\omega^2})$. We therefore obtain
$$ \rho(r,n) = e^{rT} \sup_{\omega\in \R}e^{-2n\delta\cdot x(r,\omega)} = e^{rT-2n\delta\sqrt{r}}. $$
The best choice of $r$ is the one that minimizes $\ln(\rho(r,n)) = rT - 2n\delta\sqrt{r}$, i.e., when $r = r^* = (n\delta/T)^2$. For this value, we have
$ \rho(r^*,n) = e^{-n^2\delta^2/T}. $ Thus, we deduce the estimate
\begin{equation}
\label{eq:Dirichlet_bound}
 \|E_2^n(\delta,\cdot)\|_{L^2(0,T)} \leq e^{-n^2\delta^2/T}\|E_2^0(\delta,\cdot)\|_{L^2(0,T)}. 
\end{equation}
Note that $r$ is not a parameter
of the algorithm, but an analysis parameter that we optimize to get the sharpest bounds possible for this technique.

Recall the well-known
superlinear convergence bound in \cite[Cor.~5.2]{GanderZhao}, which can be adapted to the alternating
SWR method to yield
$$ \|E_j^n\|_{L^\infty(\Omega_j\times (0,T))} \leq \mathrm{erfc}\left(\frac{n\delta}{\sqrt{T}}\right) \max_{\ell=1,2}\|E_\ell^0\|_{L^\infty(\Omega_j\times (0,T))}. $$

Since $\mathrm{erfc}(x) = e^{-x^2}\left(\frac1{x\sqrt{\pi}} + O(x^{-3})\right)$, we see that \eqref{eq:Dirichlet_bound} captures the dominant order of the superlinear convergence without the need of exact formulas for the
heat kernel nor for the inverse Laplace transform! The true power of this new analysis tool will be shown in the next section, where it allows us to obtain new results for OSWR that
remained out of reach for three decades.

\section{Convergence analysis of OSWR}
We now analyze in detail the convergence of OSWR, which is significantly more involved than for SWR with classical transmission conditions. 
We will proceed as follows:
\begin{enumerate}
\item Derive the $n$-step convergence factor $\rho(r,p,n)$, analogous to the expression \eqref{eq:classical_factor} for classical SWR, but with an extra dependence on the
Robin parameter $p$. The formula for $\rho(r,p,n)$ will contain a supremum over all frequencies $\omega$, just like in \eqref{eq:classical_factor}.
\item Calculate the supremum inside $\rho(r,p,n)$ to obtain upper bounds for each value of $r \geq 0$.
\item Study the behaviour of $r^* = r^*(p)$ that minimizes $\rho(r,p,n)$ for given values of $p$ and $n$.
\item For a given value of $n$, find the optimal $p$ that minimizes the $n$-step contraction factor $\rho(r^*(p),p,n)$. 
\end{enumerate}
The first three points will be addressed in subsections \ref{sec:3.1}--\ref{sec:3.3}, and the last point will be addressed in section \ref{sec:4}. 
\subsection{Derivation of the contraction factor}\label{sec:3.1}
When Robin transmission conditions are used, the iteration can be written in Laplace variables as
\begin{alignat*}{3}
&\text{On $\Omega_1 = (-\infty,\delta)$:} \qquad & s\hat{E}_1^n - \partial_{yy}\hat{E}_{1}^n &= 0, & \qquad (\partial_y\hat{E}_{1}^n + p\hat{E}_1^n)|_{y=\delta} &=(\partial_y\hat{E}_{2}^{n-1} + p\hat{E}_2^{n-1})|_{y=\delta},\\
&\text{On $\Omega_2 = (0,\infty)$:} \qquad & s\hat{E}_2^n - \partial_{yy}\hat{E}_{2}^n &= 0, & \qquad (-\partial_y\hat{E}_{2}^n + p\hat{E}_2^n)|_{y=0} &=(-\partial_y\hat{E}_{1}^{n} + p\hat{E}_1^{n})|_{y=0}.
\end{alignat*}
Defining $g^n(s) := (\partial_yE_{2}^n + pE_2^n)|_{y=\delta}$, a similar calculation as above leads to (cf.~\cite{Courvoisier})
$$ g^n(s) = e^{-2n\delta\sqrt{s}}\left[\frac{p-\sqrt{s}}{p+\sqrt{s}}\right]^{2n} g^0(s). $$
Thus, the $n$-step contraction factor for each $r\geq 0$ is given by
$$\rho(r,q,n) = e^{rT} \sup_{\omega \in \R} \left|e^{-2n\delta\sqrt{r+i\omega}}\left[\frac{p-\sqrt{r+i\omega}}{p+\sqrt{r+i\omega}}\right]^{2n}\right|.$$
Writing $\sqrt{r+i\omega} = x + ib$ as before, we obtain
\begin{align*}
\rho(r,q,n) = e^{rT} \sup_{\omega \in \R} e^{-2n\delta x}\left|\frac{p-x-ib}{p+x+ib}\right|^{2n} = e^{rT} \sup_{\omega \in \R}  e^{-2n\delta x}\left[\frac{(p-x)^2+b^2}{(p+x)^2+b^2}\right]^{K}.
\end{align*}
Note that we can replace $b^2$ by $x^2-r$, thanks to \eqref{eq:b}. Moreover, we have $x = \left(\frac12(r+\sqrt{r^2+\omega^2})\right)^{1/2}$, so that $x\geq \sqrt{r}$ for all $\omega\in \R$. Thus, we can rewrite $\rho(r,p,n)$ as
$$\rho(r,p,n) = \sup_{x\geq r} e^{rT-2n\delta x}\left[\frac{(p-x)^2+x^2-r}{(p+x)^2+x^2-r}\right]^{n}.$$
This motivates us to define the functions
\begin{align*}
h(x,r,p,n) &=\frac{\delta^2}{T}\left( rT - 2n\delta x +n\log\left[\frac{(p-x)^2+x^2-r}{(p+x)^2+x^2-r}\right]\right),\\
M(r,p,n) &= \sup_{x\geq \sqrt{r}} h(x,r,p,n),
\end{align*}
so that $\rho(r,p,n) = \exp(M(r,p,n)\cdot T/\delta^2)$. 
Note that by letting $\tilde{r} = \delta^2r$ and $\tik = n\delta^2/T$, we can write\footnote{We choose to write $\nu$ instead of $\tilde{n}$ for $n\delta^2/T$ in order to emphasize the
fact that this scaling parameter takes into account both the geometry and the iteration count. Indeed, we will study in Section \ref{sec:4} the behaviour of the best Robin parameter as $\nu\to 0$, which only makes sense if we fix $n$ and let $\delta\to 0$ or $T\to\infty$.}
\begin{align}\nonumber
h(x,r,p,n) &=\tilde{r} - 2{\tik}\delta \tix + {\tik}\log\left[ \frac{(p-x)^2+x^2-\tilde{r}/\delta^2}{(p+x)^2+x^2-\tilde{r}/\delta^2} \right]\\\nonumber
&=\tilde{r} - 2{\tik}\delta x +{\tik}\log\left[ \frac{(\delta p-\delta x)^2+\delta^2x^2-\tilde{r}}{(\delta p+\delta x)^2+\delta^2x^2-\tilde{r}} \right],\\
\intertext{so by letting $\tilde{x} = \delta x$ and ${\tip}=\delta p$, we get}
h(\tix,\tir,\tip,\tik) &= \tilde{r} - 2{\tik} \tilde{x} + {\tik}\log\left[ \frac{({\tip}-\tilde{x})^2+\tilde{x}^2-\tilde{r}}{({\tip}+\tilde{x})^2+\tilde{x}^2-\tilde{r}} \right].\label{def:hprx}
\end{align}
With this change of variables, the dependence on the time horizon $T$ and the overlap size $\delta$ has disappeared, since their effect is only present in the final calculation of the contraction factor $\rho$. We will therefore analyze \eqref{def:hprx} in the remainder of this section. 
To lighten the notation, we will
use the same name for the functions $h$ and $M$ for both the original and rescaled variables; the name of the arguments will make it clear which version is being used.

\emph{Remark:} In section \ref{sec:4}, we will look for the best parameter $p$ to optimize convergence. We will continue to use the rescaled variables, with the understanding 
that if the analysis yields $\tip^*$ as the optimal parameter for $\tik = n\delta^2/T$, then one should use $p^* = \tip^*/\delta$ in the implementation.

\subsection{Calculating the supremum}
In this subsection, we assume $\tir$, $\tip$ and $\tik$ to be fixed positive parameters and study the supremum
\begin{equation}
\label{eq:Mdef}
M(\tir,\tip,\tik) := \sup_{\tix\geq \sqrt{\tir}} h(\tix,\tir,\tip,\tik).
\end{equation}
To calculate this supremum, 
let us find the value of $\tix \in [\sqrt{\tir},\infty)$ that maximizes $h(\tix,\tir,\tip,\tik)$. A direct calculation gives
\begin{equation}
\label{eq:dhdx}
\frac{\partial{h}}{\partial{\tix}} = 2k \left(\frac{A\tir-B - \left(2\tix^2 - \tir - {\tip}\right)^2}{(\tip^2-\tir)^2 + 4\tix^2(\tix^2-\tir)}\right)
\end{equation}
where 
\begin{equation}
\label{def:AB}
 A = 4\tip +2 \tip^2 > 0,\qquad B = \tip^4 +{2\tip^3} -{\tip^2}. 
\end{equation}
Note that $A$ and $B$ depending only on $\tip$. If $A\tir - B \leq 0$, the function $\tix\mapsto h(\tix,\tir,\tip,\tik)$ is decreasing everywhere, so the maximum
of $h$ is achieved at $\tix = \sqrt{\tir}$. On the other hand, if $A\tir - B > 0$, then solving $\dfdx{h}{\tix}(\tix,\tir,\tip,\tik) = 0$ formally leads to the two roots
\begin{equation}\label{def:xstar}
 \tix_- = \left(\frac{\tir + {\tip} - \sqrt{A\tir-B}}{2}\right)^{1/2}, \qquad \tix_+ = \left(\frac{\tir + {\tip} + \sqrt{A\tir-B}}{2}\right)^{1/2}.
\end{equation}
\begin{figure}[t]
\begin{minipage}{0.48\textwidth}
\includegraphics[width=\textwidth,bb=0 0 225 150,clip]{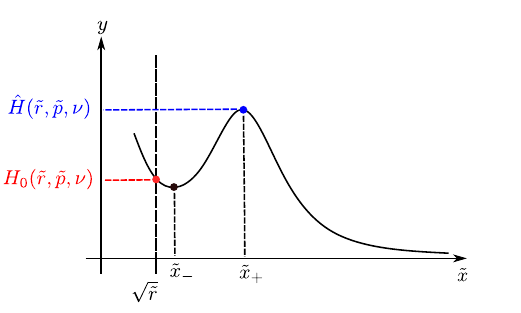}
\end{minipage}\;
\begin{minipage}{0.48\textwidth}
\includegraphics[width=\textwidth,bb=0 0 225 150,clip]{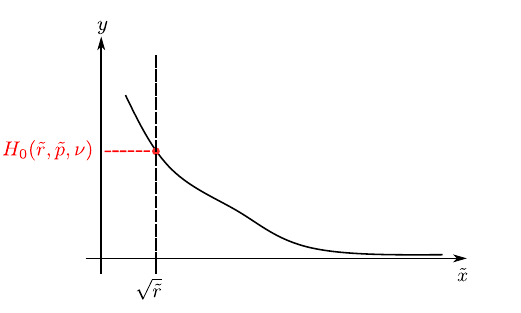}
\end{minipage}
\caption{Behaviour of $\tix\mapsto h(\tix,\tir,\tip,\tik)$ for two different values of $\tir$. Left: $A\tir - B > 0$, right: $A\tir - B < 0$. Note
that $\tix_\pm$ and $\hat{H}(\tir,\tip,\tik)$ are not defined when $A\tir-B$ is negative. \label{fig:hrx}}
\end{figure}%
We see that $\tix_+$ is always real, whereas $\tix_-$ may be real or purely imaginary. Since
$\dfdx{h}{\tix}$ is positive for $\tix_- < \tix < \tix_+$ (or simply $\tix < \tix_+$ when $\tix_-$ is purely imaginary) and negative elsewhere, $\tix\mapsto h(\tix,\tir,\tip,\tik)$ is decreasing for $\tix > \tix_+$
and has a local maximum at $\tix_+$, see Figure \ref{fig:hrx}. (The critical point $\tix_-$, when it is real, is a local minimum.) Therefore, on the interval $\tix\in [\sqrt{\tir},\infty)$, $h$ can attain its global maximum either at the end-point $\tix=\sqrt{\tir}$, or at the local maximum $\tix=\tix_+$
\emph{provided that} $\tix_+ \geq \sqrt{\tir}$, i.e., whenever
\begin{equation}
\label{r_ineq}
 \frac12\left(\tir + {\tip} + \sqrt{A\tir-B}\right) \geq \tir \iff \sqrt{A\tir-B} \geq \tir - {\tip}.
\end{equation}
Squaring both sides of \eqref{r_ineq} leads to $A\tir-B \geq (\tir-\tip)^2$, which holds when $\tir \in [\tir_-,\tir_+]$, where
$A\tir_{\pm} - B = (\tir_\pm - \tip)^2$, i.e.,
\begin{equation}
\label{def:rpm}
\tir_\pm := {\tip} + \frac{A}{2} \pm \sqrt{\frac{A^2}{4} + A\tip - B}.
\end{equation}
Note that $\tir_\pm$ are real (since $\frac14A^2-B > 0$ for all $\tip > 0$) and satisfy $\tir_+ \geq \tir_- \geq B/A$, since
\begin{equation}
\label{rmbound}
A\tir_- - B = (\tir_--\tip)^2 \geq 0 \implies \tir_- \geq B/A.
\end{equation}
Thus, the solution set of \eqref{r_ineq} contains at least the interval $[\tir_-,\tir_+]$. Additionally, if the interval $I = [B/A, \tip)$ is non-empty, then the \emph{unsquared} inequality \eqref{r_ineq} will also be satisfied for any $\tir\in I$, 
since such an $\tir$ satisfies $A\tir-B \geq 0 > \tir-\tip$. Indeed, $I$ is non-empty if and only if
\begin{equation}
\label{eq:pdelta_cond}
\frac{B}{A} < {\tip} \iff 0 < A\tip - B = {5\tip^2} - \tip^4  \iff \tip < \sqrt{5}.
\end{equation}
Thus, when $\tip < \sqrt{5}$, the solution set of \eqref{r_ineq} is $[B/A, \tip)\cup [\tir_-,\tir_+] = [B/A, \tir_+]$, since we have $B/A < \tir_- < \tip$ in this case. 
On the contrary, if $\tip \geq \sqrt{5}$, then $\tir-\tip \geq 0$ whenever $A\tir - B \geq 0$, since $\tip \leq B/A$; in this case, squaring the inequality \eqref{r_ineq} does not introduce additional solutions, so the solution set is simply $[\tir_-,\tir_+]$.

To summarize, let us define\footnote{The zero in $H_0$ suggests ``low frequency'', and the hat in $\hat{H}$ resembles a local maximum.}
$$H_0(\tir,\tip,\tik) = h(\sqrt{\tir},\tir,\tip,\tik), \qquad \hat{H}(\tir,\tip,\tik) = h(\tix_+(\tir),\tir,\tip,\tik).$$ 
It follows that for each $\tir\geq 0$, the supremum in \eqref{eq:Mdef} is given by
\begin{equation}
\label{eq:Mchar}
M(\tir,\tip,\tik) = \begin{cases}
\max (H_0(\tir,\tip,\tik), \hat{H}(\tir,\tip,\tik)), & \tir_{\min} \leq \tir \leq \tir_+,\\
H_0(\tir,\tip,\tik), & \textrm{otherwise,}
\end{cases}
\end{equation}
with
\begin{equation}
\label{def:rmin}
\tir_{\min} = \begin{cases}
\max(0, B/A), & \tip\ < \sqrt{5},\\
\tir_-, \quad & \tip \geq \sqrt{5},
\end{cases}
\end{equation}
where $\tir_\pm$ are defined in \eqref{def:rpm}. Note that there is no need to include the additional restrictions $\tir_{\min} \geq 0$ nor $\tir_{\min} \geq B/A$ in the case of $\tip \geq \sqrt{5}$,
since $\tir_-\tir_+ = B+{\tip^2} = \tip^4 + 2\tip^3 > 0$ and $\tir_+ > 0$, so we automatically have $\tir_{\min} \geq 0$; moreover,
$\tir_-$ satisfies $A\tir_--B \geq (\tir_--\tip)^2 \geq 0$, so we must have $\tir_- \geq B/A$.
Also note that when $\tip = \sqrt{5}$, we have $\tir_- = B/A > 0$, so $\tir_{\min}$ is a continuous function of ${\tip}$.

The problem of calculating the supremum is therefore reduced to comparing $H_0$ and $\hat{H}$ inside the interval $[\tir_{\min},\tir_+]$ and taking the larger of the two. The
next lemma allows us to characterize more precisely the maximum in \eqref{eq:Mchar}.
\begin{lmm}\label{lem:tangency}
Let $\tik$ and $\tip$ be positive real constants. Then the curves $y=H_0(\tir,\tip,\tik)$ and $y=\hat{H}(\tir,\tip,\tik)$ are tangent at $\tir = \tir_{+}$, and the same is true at $\tir=\tir_-$ if $\tip \geq \sqrt{5}$.
Moreover, for any $\tip>0$, we have $H_0(\tir,\tip,\tik) < \hat{H}(\tir,\tip,\tik)$ for $\tir\in (\tir_-,\tir_+)$.
\end{lmm}
\begin{proof}
We have already shown in \eqref{rmbound} that $\tir_- \geq B/A$. Since $\tir_+ \geq \tir_-$, we deduce that $A\tir_{\pm} - B \geq 0$, and $\tix_+(\tir_\pm)$ is well defined. The definition of $\tir_\pm$
implies
$\sqrt{A\tir_{\pm}-B} =| \tir_{\pm} - \tip|. $
Since $\tir_+ \geq \tip$, we must have $$\sqrt{A\tir_+-B} = \tir_+-{\tip}
\overset{\eqref{def:xstar}}{\implies} \tix_+(\tir_{+}) = \sqrt{\tir_{+}}. $$
Thus, $H_0(\tir_+,\tip,\tik) = \hat{H}(\tir_+,\tip,\tik)$. In addition, since $\tix_+(r)$ is a root of $\dfdx{h}{\tix}$ by definition, we have
\begin{align*}
\frac{\partial \hat{H}}{\partial \tir}(\tir_+,\tip,\tik)
&= \frac{\partial h}{\partial \tir}(\tix_+(\tir_+),\tir_+,\tip,\tik) + \cancelto{0}{\frac{\partial h}{\partial \tix}(\tix_+(\tir_+),\tir_+,\tip,\tik)}\cdot \frac{d\tix_+}{d\tir}(\tir_+)\\
&= \frac{\partial h}{\partial \tir}(\sqrt{\tir_+},\tir_+,\tip,\tik) 
+ \cancelto{0}{\frac{\partial h}{\partial \tix}(\sqrt{\tir_+},\tir_+,\tip,\tik)}\cdot \frac{1}{2\sqrt{\tir_+}} = \frac{\partial H_0}{\partial \tir}(\tir_+,\tip,\tik). 
\end{align*}
Therefore, the curves $y=H_0(\tip,\tir,\tik)$ and $y=\hat{H}(\tip,\tir,\tik)$ are in fact tangent at $r=\tir_+$.

If $\tip \geq \sqrt{5}$, then $A\tip - B \leq 0$ by \eqref{eq:pdelta_cond}, which  implies $\tir_-\geq \tip$ by \eqref{def:rpm}. Thus, we again have  $\tix_+(\tir_-) = \sqrt{\tir_-}$, so tangency also holds at $\tir=\tir_-$.

Finally, we show that $\hat{H}(\tir,\tip,\tik) > H_0(\tir,\tip,\tik)$ for $\tir\in (\tir_-,\tir_+)$, regardless of the value of ${\tip}$. First, substituting
$\tix=\sqrt{\tir}$ into \eqref{eq:dhdx} leads to
$$ \frac{\partial h}{\partial \tix}(\sqrt{\tir},\tir,\tip,\tik) = 2\tik\left(\frac{A\tir-B-(\tir-{\tip})^2}{(\tip^2-\tir)^2}\right).$$
By the definition of $\tir_\pm$, we have $A\tir-B > (\tir-\tip)^2$ whenever $\tir \in (\tir_-,\tir_+)$.
Hence, $\tix\mapsto h(\tix,\tir,\tip,\tik)$ is strictly increasing in a neighbourhood of $\tix=\sqrt{\tir}$, meaning that $h(\cdot,\tir,\tip,\tik)$ cannot achieve its maximum at $\tix=\sqrt{\tir}$. Therefore, the maximum must occur at $\tix_+(r)$, so we must have $\hat{H}(\tir,\tip,\tik) > H_0(\tir,\tip,\tik)$.
\end{proof}

Using the definitions of $A$ and $B$ in  \eqref{def:AB}, one can verify by direct calculation that $\tir_{-} \leq \tip^2 \leq \tir_{+}$. The above lemma implies that when $\tir= \tip^2$, the supremum is always given by
$\hat{H}(\tir,\tip,\tik)$, so the fact that $H_0(\tir,\tip,\tik) = -\infty$ when $\tir=\tip^2$ poses no difficulty. It also implies that when $\tip \geq \sqrt{5}$, the maximum in \eqref{eq:Mchar} is always
given by $\hat{H}$. We now treat the case of $\tip < \sqrt{5}$.
\begin{lmm}
Let $0 < \tip < \sqrt{5}$. Then $H_0(\tir_{\min},\tip,\tik) \geq \hat{H}(\tir_{\min},\tip,\tik)$. Moreover, the function $\tir\mapsto \hat{H}(\tir,\tip,\tik) - H_0(\tir,\tip,\tik)$ is strictly increasing on $\tir\in [\tir_{\min}, \tir_-]$, and there is a unique crossing point $\tir^\times \in [\tir_{\min},\tir_-]$ at which $H_0(\tir^\times,\tip,\tik) = \hat{H}(\tir^\times,\tip,\tik)$.
\end{lmm}
\begin{proof}
We have either $\tir_{\min} = B/A$ or $\tir_{\min} = 0$, depending on the sign of $B$. If $B > 0$, then $\tir_{\min} = B/A$, and we have
$$ \frac{\partial h}{\partial \tix}(\sqrt{\tir_{\min}},\tir_{\min},\tip,\tik) = 2\tik\left(\frac{0-(\tir_{\min}-{\tip})^2}{(\tip^2-\tir_{\min})^2}\right) < 0,$$
so that $\tix \mapsto h(\tix,\tir_{\min},\tip,\tik)$ is decreasing over all $\tix\in [\tir_{\min},\infty)$. Therefore, we must have $\hat{H}(\tir_{\min},\tip,\tik) = h(\tix_+(\tir_{\min}),\tir_{\min},\tip,\tik)) \leq
h(\sqrt{\tir_{\min}},\tir_{\min},\tip,\tik) = H_0(\tir_{\min},\tip,\tik)$. If $\tir_{\min} = 0$, i.e., if $B \leq 0$, then a direct calculation shows for all $\tix \geq 0$ that
$$ h(\tix,0,\tip,\tik) = -2\tik\tix + \tik\log\left(\frac{(\tip-\tix)^2+\tix^2}{(\tip+\tix)^2+\tix^2}\right) \leq 0 = h(0,0,\tip,\tik),$$
since the argument of the logarithm is less than or equal to 1 for any positive $\tip$ and non-negative $\tix$. Letting $\tix = \tix_+(0)$ shows that $\hat{H}(0,\tip,\tik) \leq H_0(0,\tip,\tik)$, as required.
Finally, Lemma \ref{lem:tangency} shows that $\hat{H}(\tir_-,\tip,\tik) \geq H_0(\tir_-,\tip,\tik)$; since both $\tir\mapsto \hat{H}(\tir,\tip,\tik)$ and $\tir\mapsto H_0(\tir,\tip,\tik)$ are continuous functions,
the intermediate value theorem implies the curves must cross at least once at $\tir^\times$, somewhere between $\tir_{\min}$ and $\tir_-$. 
\begin{figure}[t]
\includegraphics[width=0.48\textwidth,bb=110 235 490 545,clip]{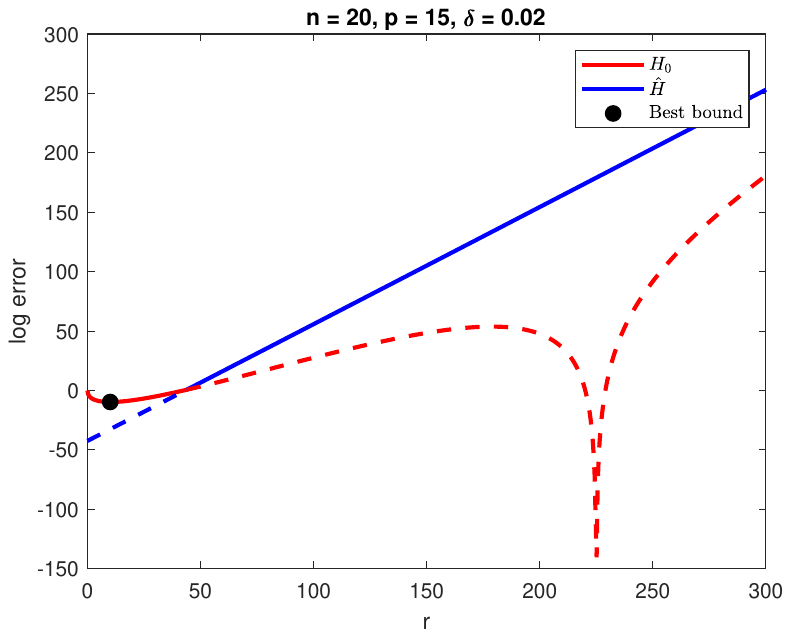}
\includegraphics[width=0.48\textwidth,bb=110 235 490 545,clip]{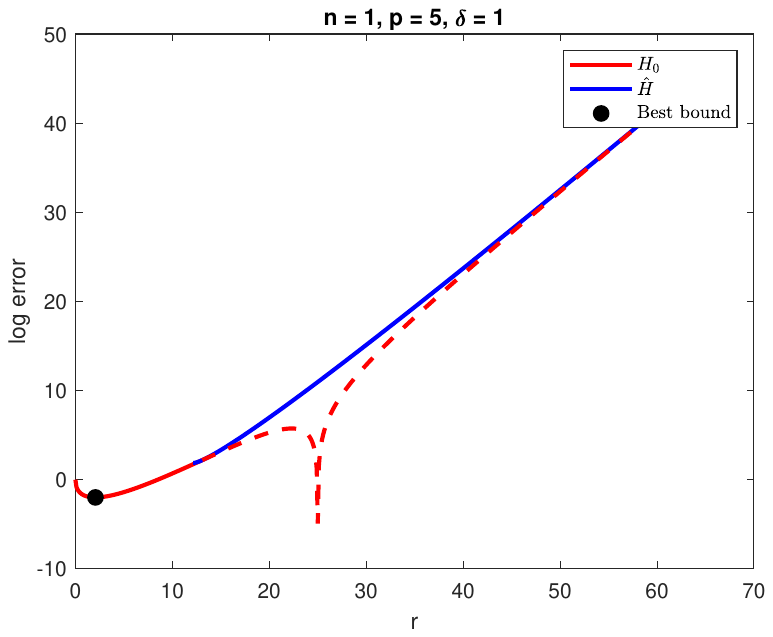}\\
\includegraphics[width=0.48\textwidth,bb=110 235 490 545,clip]{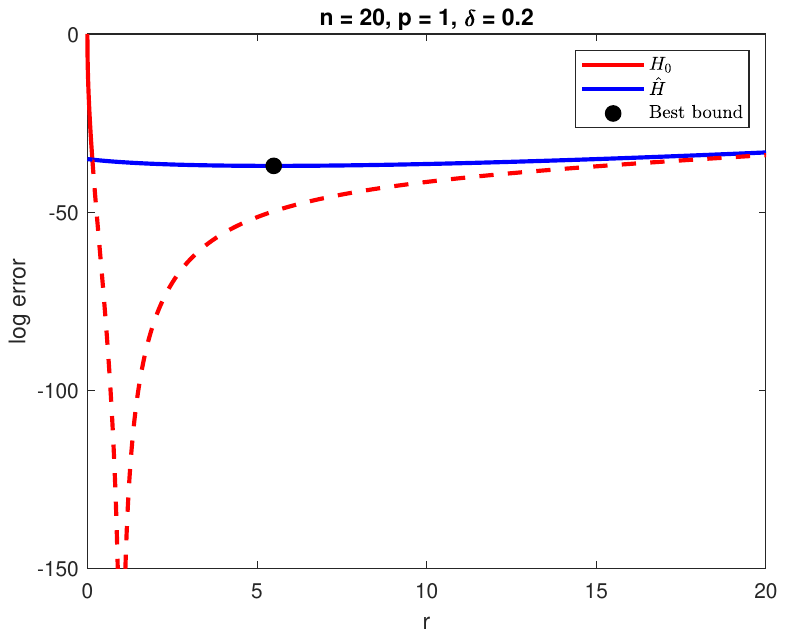}
\includegraphics[width=0.48\textwidth,bb=110 235 490 545,clip]{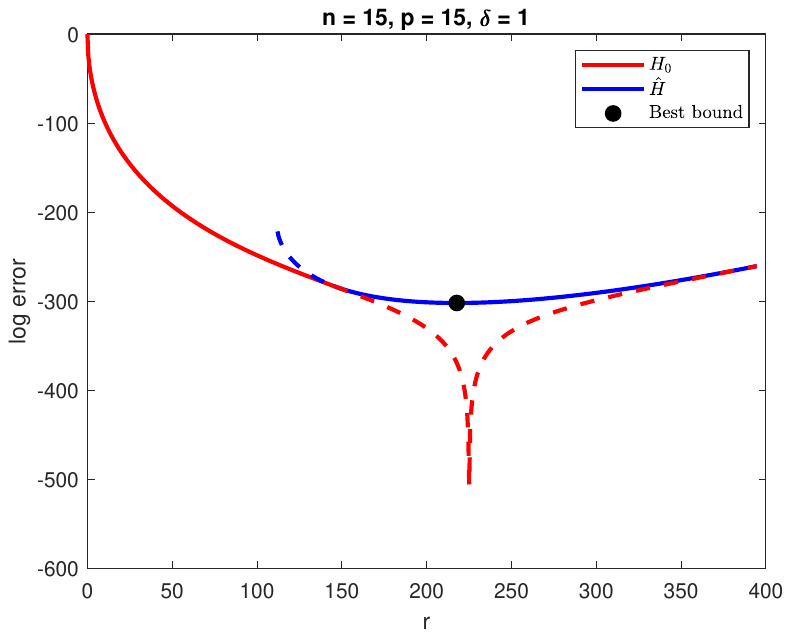}\\
\includegraphics[width=0.48\textwidth,bb=110 235 490 545,clip]{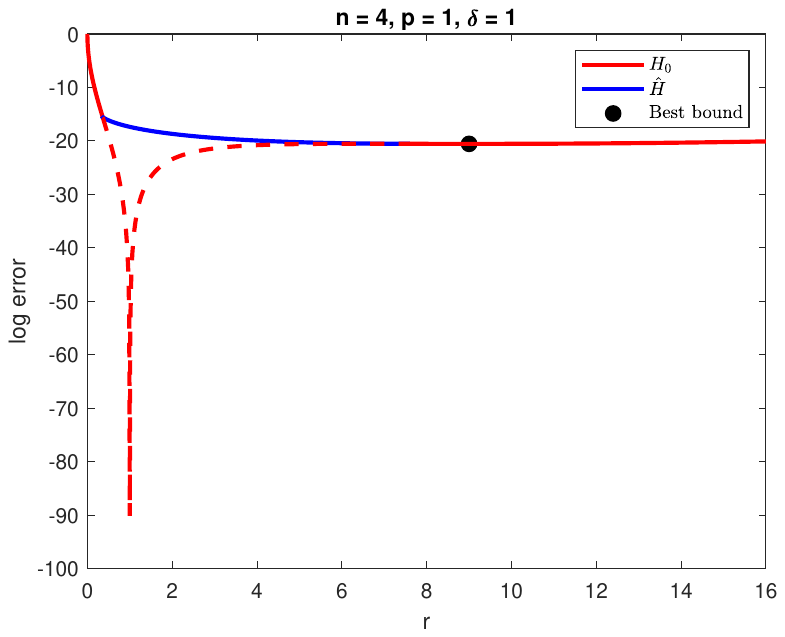}
\includegraphics[width=0.48\textwidth,bb=110 235 490 545,clip]{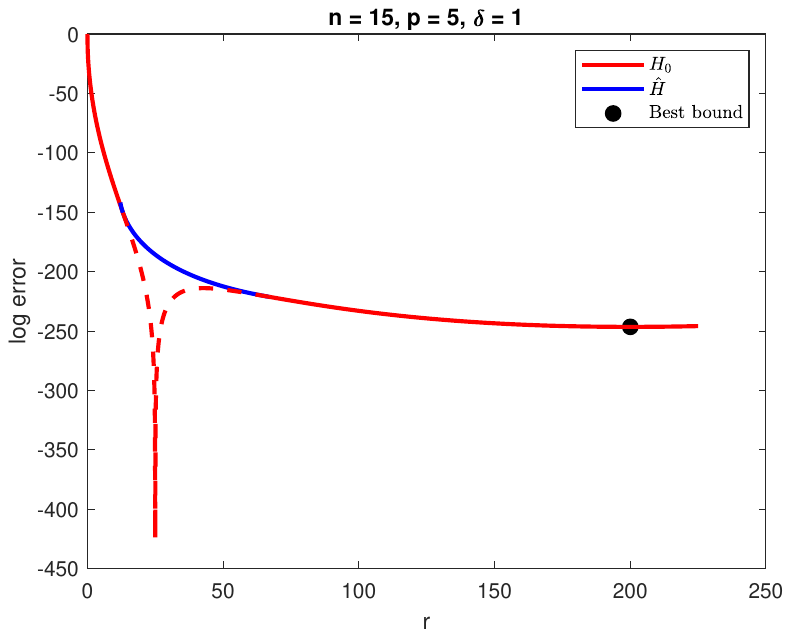}\\
\includegraphics[width=0.48\textwidth,bb=110 235 490 545,clip]{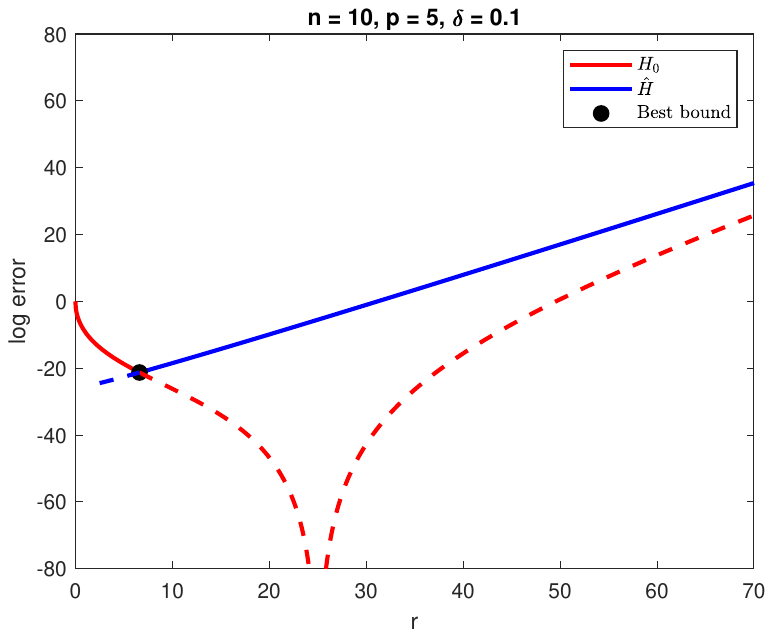}
\includegraphics[width=0.48\textwidth,bb=110 235 490 545,clip]{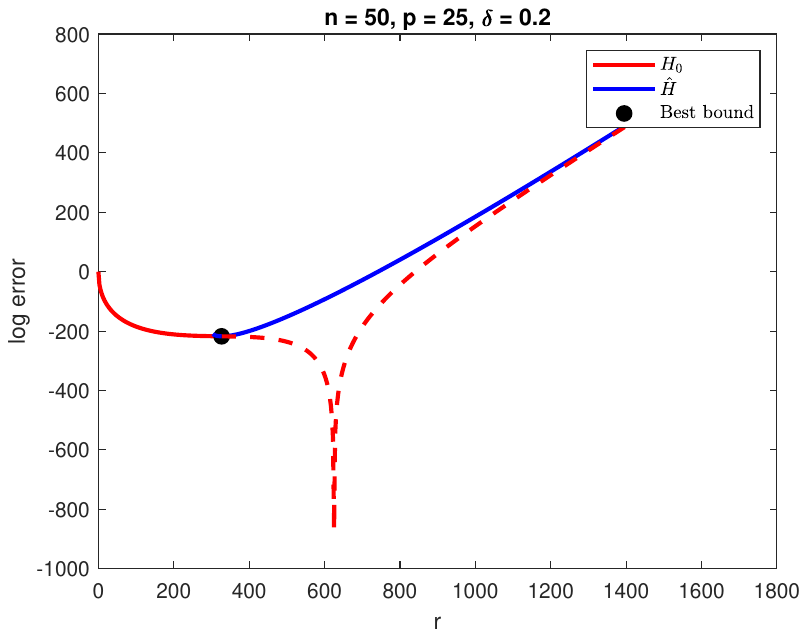}
\caption{Different configurations for the $r$-dependent error bound. The bound
is defined by the union of the solid curves, and the dashed curves represent
portions of $H_0$ or $\hat{H}$ that do not contribute to the error bound.
The black dot indicates the global minimum of $r\mapsto M(r,p,n)$. The left column
illustrates cases where ${\tip} = p\delta < \sqrt{5}$, while the right column corresponds to ${\tip} = p\delta \geq \sqrt{5}$.
Only the bottom left panel is characterized by equioscillation; the other minima are all critical points of $H_0$ or $\hat{H}$ (or both,
as is the case for the bottom right panel).\label{fig4.1}}
\end{figure}
\afterpage{\clearpage}

We now show that $\tir\mapsto \hat{H}(\tir,\tip,\tik) - H_0(\tir,\tip,\tik)$ is strictly increasing on $\tir\in [\tir_{\min},\tir_-]$, which implies that the function cannot vanish for more than one value of $\tir$. A direct calculation shows that
\begin{equation}
\label{lem4.2b}
\dfdx{\hat{H}}{\tir}(\tir,\tip,\tik) - \dfdx{H_0}{\tir}(\tir,\tip,\tik) = \frac{\tik}{\sqrt{\tir}} + \frac{2\tik\tip}{\sqrt{\tir}(\tip^2-\tir)} - \frac{4\tik\tip\tix_+}{G},
\end{equation}
where $\tix_+ = \tix_+(\tir)$ and $G = G(\tix_+,\tir,\tip) = (\tip^2-\tir)^2+4\tix_+^2(\tix_+^2-\tir) > 0$, since $\tix_+ \geq \sqrt{\tir}$. 
Since $\tix_+$ is a local maximum of $\tix\mapsto h(\tip,\tir,\tix)$, we have
$$ 0 = \dfdx{h}{\tix}(\tix_+,\tir,\tip,\tik) = -2\tik - \frac{4\tik\tip(\tip^2-2\tix_+^2-\tir)}{G} \implies G = {2\tip}(2\tix_+^2+\tir-\tip^2). $$
Substituting this into \eqref{lem4.2b} leads to
\begin{align*}
\dfdx{\hat{H}}{\tir}(\tir,\tip,\tik) - \dfdx{H_0}{\tir}(\tir,\tip,\tik) &= \frac{2\tik\tip}{G\sqrt{\tir}}\left(2\tix_+^2+\tir-\tip^2 + \frac{G}{\tip^2-\tir} - 2\tix_+\sqrt{\tir}\right)\\
&= \frac{2\tik\tip}{G\sqrt{\tir}(\tip^2-\tir)}\left(G - (\tip^2-\tir)^2 + (2\tix_+^2 - 2\tix_+\sqrt{\tir})(\tip^2-\tir)\right)\\
&=  \frac{2\tik\tip}{G\sqrt{\tir}(\tip^2-\tir)}\left(4\tix_+^2(\tix_+^2-\tir) + 2\tix_+(\tix_+-\sqrt{\tir})(\tip^2-\tir)\right) > 0,
\end{align*}
since $\tir \leq \tir_- < \tip^2$. Thus, $\tir\mapsto \hat{H}(\tip,\tir)-H_0(\tip,\tir)$ is strictly increasing, as required.
\end{proof}
The following theorem summarizes the behaviour of $M(\tir,\tip,\tik)$ for all values of $\tip$.
\begin{thrm}\label{thm:M_prop}
Let $\tik > 0$ and $\tip >0$ be given. Let $M(\tir,\tip,\tik)$ be the bound on the convergence factor defined in \eqref{eq:Mdef}. Then for $\tip \geq \sqrt{5}$, we have
$$ M(\tir,\tip,\tik) = \begin{cases}
\hat{H}(\tir,\tip,\tik) = h(\tix_+(\tir),\tir,\tip,\tik), \quad & \tir_- \leq \tir \leq \tir_+,\\
H_0(\tir,\tip,\tik) = h(\sqrt{\tir},\tir,\tip,\tik), & \text{otherwise},
\end{cases}$$
where $\tir_-$ and $\tir_+$ are defined in \eqref{def:rpm}. Moreover, $\tir\mapsto M(\tip,\tir)$ is continuously differentiable on $\tir\in [0,\infty)$. If $\tip < \sqrt{5}$, then
$$ M(\tir,\tip,\tik) = \begin{cases}
\hat{H}(\tir,\tip,\tik) = h(\tix_+(\tir),\tir,\tip,\tik), \quad & \tir^\times \leq \tir \leq \tir_+,\\
H_0(\tir,\tip,\tik) = h(\sqrt{\tir},\tir,\tip,\tik), & \text{otherwise},
\end{cases}$$
where $\tir^\times$ is the unique point in $[\tir_{\min},\tir_-]$ that satisfies  $\hat{H}(\tir^\times,\tip,\tik) = H_0(\tir^\times,\tip,\tik)$. Moreover, $\tir\mapsto M(\tir,\tip,\tik)$ is continuous on $\tir\in [0,\infty)$ and continuously differentiable everywhere except at $\tir=\tir^\times$. 
\end{thrm}
\subsection{Minimizing the bound for fixed $\tip$ and $\tik$}\label{sec:3.3}
We now turn our attention to finding the parameter $\tir$ that yields the sharpest bound, i.e., the one that minimizes $M(\tir,\tip,\tik)$. To show the existence of such an $\tir$, we observe that $\tir\mapsto M(\tir,\tip,\tik)$ is increasing for large enough $\tir$. Indeed, since
\begin{equation}
\label{eq:dH0dr}
\frac{\partial H_0}{\partial \tir} = \frac{\partial}{\partial \tir} \left[\tir - 2\tik\sqrt{\tir} + \tik\log\left(\frac{(\tip-\sqrt{\tir})^2}{(\tip+\sqrt{\tir})^2}\right)\right] = 1 - \frac{\tik}{\sqrt{\tir}} + \frac{2\tik\tip}{\sqrt{\tir}(\tir-\tip^2)}, 
\end{equation}
we have $\dfdx{M}{\tir} > 0$ for all $\tir \geq \max(\tik^2,\tir_+)$. Moreover, $\tir\mapsto H_0(\tir,\tip,\tik)$ is decreasing in a neighbourhood containing $\tir=0$ (in fact, $\dfdx{M}{\tir} \to -\infty$ as $\tir\to 0^+$). Therefore, $\tir\mapsto M(\tip,\tir)$ has a global minimum somewhere on the interval $(0, \max(\tik^2,\tir_+))$.
If ${\tip} \geq 5$, the minimum must be either a critical point of $H_0$ or $\hat{H}$, since $M$ is continuously differentiable. If ${\tip} < 5$, the minimum can be a critical point, or it can be the intersection point $\tir^\times \in [\tir_{\min}, \tir_-]$, where $H_0(\tir^\times,\tip,\tik) = \hat{H}(\tir^\times,\tip,\tik)$. At such a point, we have \emph{equioscillation} between the contraction
factor at $\tix=\sqrt{\tir}$ and that at the local maximum $\tix=\tix_+(\tir)$. Figure \ref{fig4.1} illustrates the different possible locations for the minimum, depending on $\tip$ and $\tik$.
In all but one case, the minimum corresponds to a critical point of $H_0$ or $\hat{H}$ (or both, when they coincide); the remaining case is characterized by equioscillation, i.e., we have $\tir=\rx$. This leads us to
study the critical points of $H_0$ and $\hat{H}$.
To study the critical points of $H_0$, we let $z = \sqrt{\tir}$ and rewrite \eqref{eq:dH0dr} as follows:
$$\dfdx{H_0}{\tir} = \frac{\tik}{\tip z}\left(\frac{\tip z}{\tik} - \left[{\tip} + \frac{2\tip^2}{\tip^2-z^2}\right]\right) =: \frac{\tik}{\tip z}\left(\frac{\tip z}{\tik} - F_{\tip}(z)\right), $$
where we observe that $F_{\tip}$ is independent of $\tik$. We now have the following lemma.
\begin{figure}
{
\centering
\begin{overpic}[width=0.6\textwidth,bb=120 265 480 535,clip]{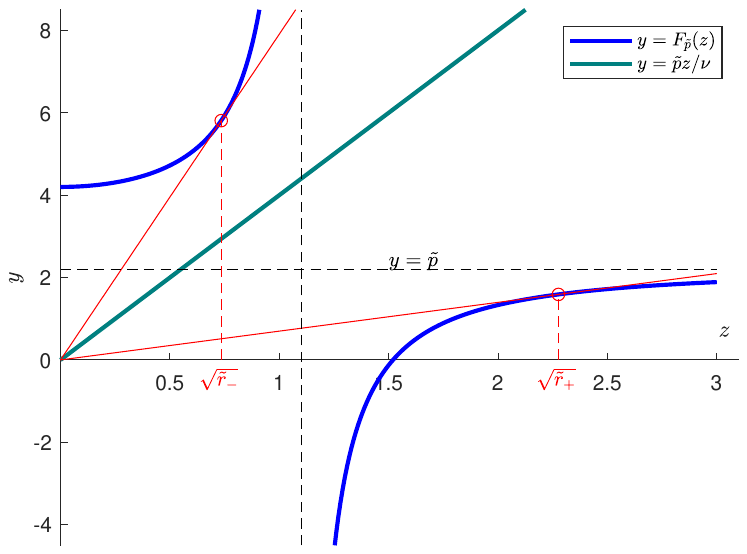}
\put(60,66){\Huge$\lcurvearrowse$}
\put(62,61){\small Increasing $\tik$}
\end{overpic}
\caption{Plot of $y=F_{\tip}(z)$ and $y=\tip z/\tik$ in the case where $\tip/F_{\tip}'(\sqrt{\tir_-}) < \tik < \tip/F_{\tip}'(\sqrt{\tir_+})$.\label{fig:crit_H0}}}
\end{figure}%
\begin{lmm}[Critical points of $H_0$]\label{lem:H0crit}
Let $\tir_-$ and $\tir_+$ be defined as in \eqref{def:rpm} and $\tik > 0$.
\begin{itemize}
\item If $\tik < \tip/F_{\tip}'(\sqrt{\tir_-})$, then $H_0$ has a local minimum $\tir_*$ and a local maximum $\tir^*$ with $\tir_* < \tir_- < \tir^* < \tip^2$ and no other critical points in $[0,\infty)$.
\item If $\tik > \tip/F_{\tip}'(\sqrt{\tir_+})$, then $H_0$ has a local minimum $\tir_*$ and a local maximum $\tir^*$ with $\tip^2 < \tir^*<\tir_+<\tir_*$ and no other critical points in $[0,\infty)$.
\item If $\tik=\tip/F'_{\tip}(\sqrt{\tir_-})$ (resp.~$\tik=\tip/F_{\tip}'(\sqrt{\tir_+})$), then $H_0$ has a unique critical point that is also an inflection point at $\tir_-$ (resp. $\tir_+$).
\item If $\tip/F_{\tip}'(\sqrt{\tir_-}) < \tik < \tip/F_{\tip}'(\sqrt{\tir_+})$, then $H_0$ has no critical points in $[0,\infty)$.
\end{itemize}
\end{lmm}
\begin{proof}
We start by rewriting $\tir_\pm$. Substituting the definitions of $A$ and $B$ in \eqref{def:AB} into \eqref{def:rpm} leads to
\begin{equation}
\label{eq:rpm_pd}
\begin{aligned}
\tir_\pm &= {3\tip} + \tip^2  \pm\sqrt{\left(\tip^4 + {4\tip^3} + {4\tip^2}\right) + {4\tip^2} + {2\tip^3} - \left(\tip^4 + {2\tip^3} -
{\tip^2}\right)}\\
&= {\tip}\left(3 + {\tip} \pm \sqrt{9 + 4{\tip}}\right).
\end{aligned}
\end{equation}
Next, we consider the zeros of $\dfdx{H_0}{\tir}$, which occur when there are intersections between $y = F_{\tip}(z)$ and the straight line $y = \tip z/\tik$, see Figure \ref{fig:crit_H0}. 
As $\tik$ varies from $0$ to $\infty$, we see that the straight line intersects the upper branch of $y=F_{\tip}(z)$ if $\tik$ is small; it intersects the lower branch if $\tik$ is large, and there are no intersections for intermediate values of $\tik$. The transition between these three regimes occurs when the straight line is tangent to one of the branches, see the red lines in Figure  \ref{fig:crit_H0}. Therefore, it suffices to show that these
points of tangency occur at $z=\sqrt{\tir_\pm}$. Indeed, for a given $z_0 \neq \tip$, the tangent line is of the form
$$ y = F_{\tip}(z_0) + F_{\tip}'(z_0)(z-z_0) = {\tip} + \frac{2\tip^2}{\tip^2-z_0^2} + \frac{4\tip^2z_0}{(\tip^2-z_0^2)^2}(z-z_0).$$
This straight line passes through the origin if
$$ 0 =  {\tip} + \frac{2\tip^2}{\tip^2-z_0^2} - \frac{4\tip^2z_0^2}{(\tip^2-z_0^2)^2} = \frac{{\tip}(\tip^2-z_0^2)^2 + 2\tip^2(\tip^2-z_0^2) -4\tip^2z_0^2}{(\tip^2-z_0^2)^2}. $$
Thus, we see that the numerator is zero whenever
$$ z_0^4 -2z_0^2(\tip^2 + 3\tip) + (\tip^4 + 2\tip^3) = 0. $$
Solving this quadratic equation in $z_0^2$ leads to
$$ z_0^2 = \tip^2 + {3\tip} \pm \sqrt{\left(\tip^4+{6\tip^3} + {9\tip^2}\right) - \left(\tip^4 + {2\tip^3}\right)} = {\tip}\left(3 + {\tip}\pm
\sqrt{9 + 4{\tip}}\right) = \tir_{\pm}.$$
The critical values of $\tik$ are therefore given by $\tip/\tik = F_{\tip}'(z_0)$ with $z_0 = \sqrt{\tir_{\pm}}$, and the result follows.
\end{proof}
We now turn our attention to the function $\tir\mapsto \hat{H}(\tir,\tip,\tik)$.
\begin{lmm}[Critical points of $\hat{H}$]
For each fixed $\tip > 0$, $\tir \mapsto \hat{H}(\tir,\tip,\tik)$ has at most one critical point $\tir^* \in [\tir_{\min},\tir_+]$ such that $\dfdx{\hat{H}}{\tir}(\tir^*,\tip,\tik) = 0$.
\end{lmm}
\begin{proof}
If $\tir^*$ is a critical point of $\hat{H}(\tir,\tip,\tik) = h(\tix_+(\tir),\tir,\tip,\tik)$, 
then
\begin{equation}
\label{eq:dHHdr}
\dfdx{\hat{H}}{\tir}(\tir^*,\tip,\tik) = \dfdx{h}{\tir}(\tix_+(\tir^*),\tir^*,\tip,\tik) + \underbrace{\dfdx{h}{\tix}(\tix_+(\tir^*),\tir^*,\tip,\tik)}_{=0}\frac{d\tix_+}{d\tir}(\tir^*) = 0.
\end{equation}
Therefore, $(\tix_+,\tir^*)$ is a solution of the system
\begin{equation}
\label{lem_eq0}
\left\{
\begin{aligned}
\dfdx{h}{\tir} &=& 1-\frac{4\tik\tip\tix}{G(\tix,\tir,\tip)}\qquad\qquad &= 0,\\
\dfdx{h}{\tix} &=& -2\tik - \frac{4\tik\tip(\tip^2-2\tix^2-\tir)}{G(\tix,\tir,\tip)} &= 0,
\end{aligned}
\right.
\end{equation}
where $G(\tix,\tir,\tip) = (\tip^2-\tir)^2 + 4\tix^2(\tix^2-\tir) > 0$ for all $\tix > \sqrt{\tir}$. At the solution $(\tix_+,\tir^*)$, we can isolate $G(\tix_+,\tir^*,\tip)$ from each equation to get
\begin{equation}\label{eq:r_stationary}
 G(\tix_+,\tir^*,\tip) = 4\tik\tip\tix_+ = -{2\tip}(\tip^2-2\tix_+^2-\tir^*), 
\end{equation}
from which we obtain
\begin{equation}\label{eq:r_from_x}
\tir^* = \tip^2 + 2\tik \tix_+ - 2\tix_+^2. 
\end{equation}
Hence, we define $\tir = \varphi(\tix) = \tip^2 + 2\tik\tix - 2\tix^2.$ Next, we use the explicit formula of $\tix_+$ in \eqref{def:xstar} to define the function $\tir = \psi(\tix)$, such that
$$ \tix_+ = \left(\frac{\psi(\tix_+) + {\tip} + \sqrt{A\psi(\tix_+) - B}}{2}\right)^{1/2} =: X(\psi(\tix_+)). $$
In other words, $\psi$ and $X$ are inverse functions. We claim that $\psi$ is an increasing convex function. Indeed, since $X(\tir) > 0$ and
$$
2(X(\tir))^2 = \tir + {\tip} + \sqrt{A\tir-B} \implies 4X(\tir)X'(\tir) = 1 + \frac{A}{2\sqrt{A\tir-B}} > 0 \implies X'(\tir) > 0, 
$$
Differentiating once more gives
$$ 4(X'(\tir))^2 + X(\tir)X''(\tir) = \frac{-A^2}{4(A\tir-B)^{3/2}} \implies X''(\tir) = \frac{1}{X(\tir)}\left(\frac{-A^2}{4(A\tir-B)^{3/2}} - 4(X'(\tir))^2\right) < 0. $$
Since the inverse of an increasing concave function is increasing and convex, we conclude that $\psi(\tix)$ is increasing and convex.
We now show that there can be at most one intersection between $\tir = \varphi(\tix)$ and $\tir=\psi(\tix)$ within the interval $\tix \in [\tix_{\min},x^+] := [X(\tir_{\min}),X(\tir^+)].$
First, we either have $\tir_{\min} = 0$ or 
$$ \tir_{\min} = \frac{B}{A} = \frac{\tip^4 + {2\tip^3} - {\tip^2}}{2\tip^2+{4\tip}} < \frac{\tip^4+2\tip^3}{2\tip^2+{4\tip}} = \frac{\tip^2}{2}.$$
So either way, $\tir_{\min} < \tip^2 = \varphi(0)$. If $\tir_{\min} > \varphi(\tix_{\min})$, then 
$$ \tip^2 > \tip^2 + 2\tik \tix_{\min} - 2\tix_{\min}^2 \implies \tix_{\min} > \tik. $$
Since $\varphi'(\tix) = 2\tik - 4\tix < 0$ for $\tix > \tik$, we see that $\psi - \varphi$ is strictly increasing on $[\tix_{\min},\tix_+]$, so any $\tix$ in this interval must satisfy
$$\psi(\tix) - \varphi(\tix) \geq \psi(\tix_{\min}) - \varphi(\tix_{\min}) = \tir_{\min} - 
\varphi(\tix_{\min}) > 0.$$ In other words, no intersection can exist between $\tir=\psi(\tix)$ and $\tir = \varphi(\tix)$. Now suppose $\tir_{\min} \leq \varphi(\tix_{\min})$. Since $\psi(\tix_{\min})-\varphi(\tix_{\min}) \leq 0$ and $\psi - \varphi$ is convex on $[\tix_{\min},\tix_+]$, it can have at most one root $\tix_+ \in [\tix_{\min},\tix_+]$. 
\par In summary, since there is at most one intersection between $\tir=\varphi(\tix)$ and $\tir=\psi(\tix)$, there is at most one solution to the optimality system \eqref{lem_eq0}, so $\hat{H}$
can have at most one critical point.
\end{proof}

\begin{figure}[t]
\begin{center}
\includegraphics[width=0.48\textwidth,bb=110 235 490 545,clip]{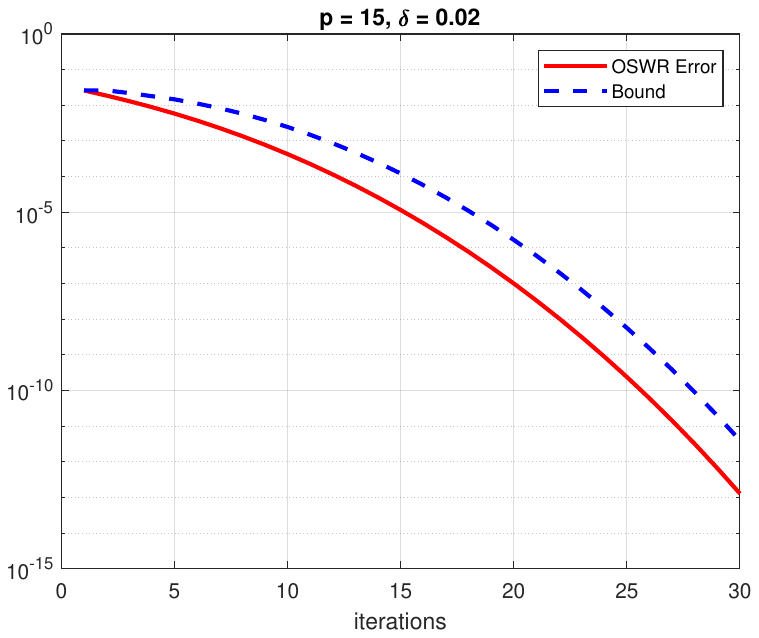}
\includegraphics[width=0.48\textwidth,bb=110 235 490 545,clip]{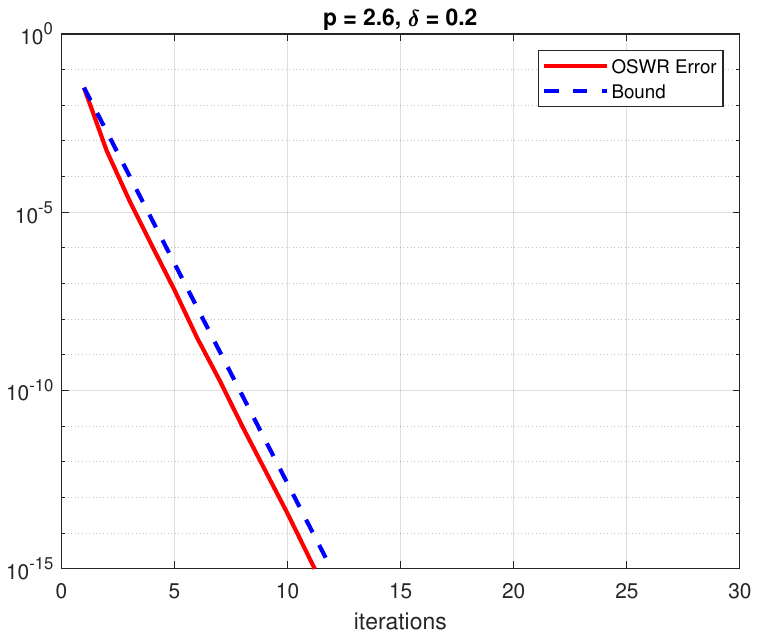}\\
\includegraphics[width=0.48\textwidth,bb=110 235 490 545,clip]{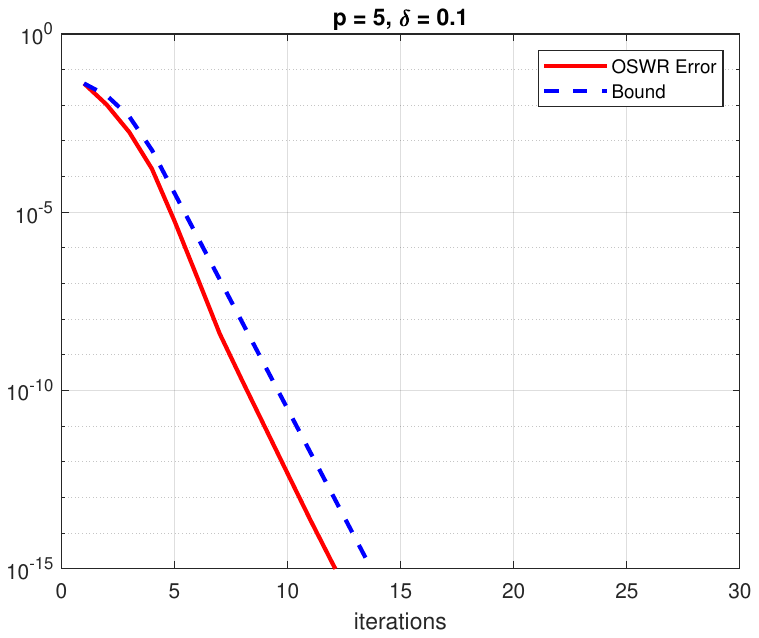}
\includegraphics[width=0.48\textwidth,bb=110 235 490 545,clip]{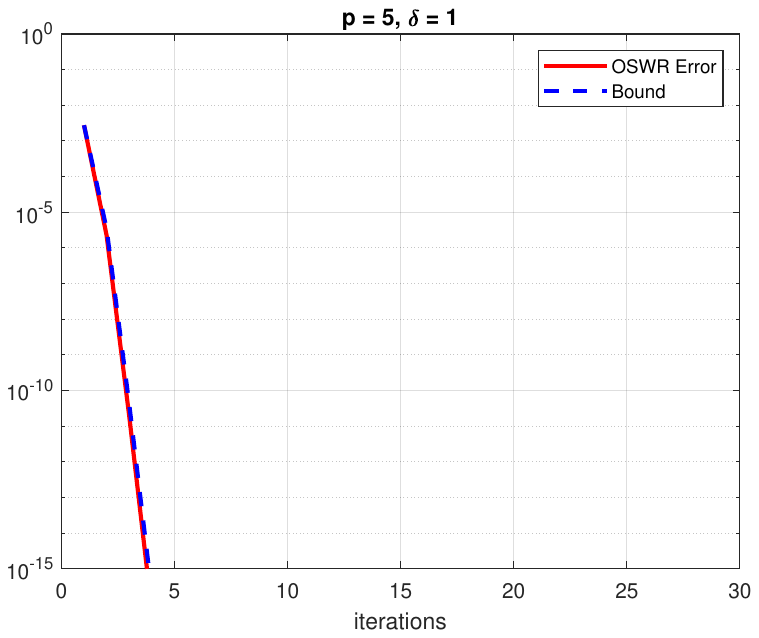}
\caption{Comparison between convergence curves for OSWR and the superlinear bound for different values of $p$ and $\delta$. We have $T=1$ in all cases. \label{fig3}}
\end{center}
\end{figure}
Given the different ways in which $\tir\mapsto M(\tir,\tip,\tik)$ can achieve its minimum (see Figure \ref{fig4.1}), it is difficult to obtain a closed form expression for the optimal $\tir$ and
the resulting convergence bound analytically. Nonetheless, the optimal bound is
easy to calculate numerically for a given $\tip = \delta p$:
\begin{enumerate}
\item If ${\tip} \geq \sqrt{5}$:
\begin{itemize}
\item Compute $\tir_\pm$, defined in \eqref{def:rpm}. If $H_0$ has a local minimum $\tir_0$ by Lemma \ref{lem:H0crit}, then this is also the global minimum of $\tir\mapsto M(\tir,\tip,\tik)$.
\item If $H_0$ has no local minimum, then find the local minimum $\hat{r}$ of $\hat{H}$ by solving \eqref{def:xstar} and \eqref{eq:r_from_x} simultaneously.
\end{itemize}
\item If ${\tip} < \sqrt{5}$: 
\begin{itemize}
\item Compute the local minimum $\widehat{r}$ of $\hat{H}$ by solving \eqref{def:xstar} and \eqref{eq:r_from_x} simultaneously. If $\widehat{r}$ exists and satisfies $H_0(\widehat{r},\tip,\tik) \leq \hat{H}(\widehat{r},\tip,\tik)$, then $\widehat{r}$ is the global minimum of $\tir\mapsto M(\tir,\tip,\tik)$.
\item Otherwise, compute the equioscillation point $\rx$ by setting $H_0$ equal to $\hat{H}$ (e.g.~using $\tir_-$ as the starting point for Newton). If $\dfdx{H_0}{\tir}\dfdx{\hat{H}}{\tir} < 0$ at $\tir=\rx$, then $\rx$ is the global minimum of $M$.
\item If neither $\widehat{r}$ nor $\rx$ works, then $H_0$ must have a local minimum $\tir_0$ that coincides with that of $M$. Calculate $\tir_0$ by finding
the root of $\dfdx{H_0}{\tir} = 0$ that is outside the interval $[\tir_-,\tir_+]$, see \eqref{eq:dH0dr} for an expression of $\dfdx{H_0}{\tir}$.
\end{itemize}
\end{enumerate}
We use this procedure to compute bounds for the four test cases shown in Figure \ref{fig3}. The convergence curves correspond to OSWR applied to the 1D homogeneous heat equation
posed on the \emph{finite} domain $y \in \Omega = (-4,4)$ and integrated on the time interval $t \in [0,1]$. We impose homogeneous Dirichlet conditions at $y=\pm 4$ and zero initial conditions.
The problem is discretized using backward Euler in time and centered finite differences in space, with $\Delta y = 1/800$ and $\Delta t = 1/2400$. The two subdomains are of
equal size, with the overlap indicated on the plot title; the Robin parameters used in the transmission conditions are also shown in each case. Random initial guesses are used for
the interface at iteration 0. We see that in all cases, our theoretical bounds track the convergence curves closely and are able to capture the superlinear behaviour reliably.
\afterpage{\clearpage}

\section{Optimizing the Robin parameter}\label{sec:4}
Interestingly, the analysis becomes more tractable when we seek an optimal Robin parameter
$\tip$, since optimality yields an extra equation that can be used to obtain asymptotic 
estimates under various scalings of $\tik$. We start by showing that
such an optimal $\tip$ exists.

\begin{thrm}
Let $\tik > 0$ be given. For any $\tip> 0$, let $\tir_* = \tir_*(\tip,\tik)$ be a global minimum of the function $\tir\mapsto M(\tir,\tip,\tik)$. Then there exists a value $\tipm$,
which depends on $\tik$, such that for any $\tip > \tipm$, we have
$M(\tir_*(\tip,\tik),\tip,\tik) > M(\tir_*(\tipm,\tik),\tipm,\tik)$.
As a result, the function $\tip\mapsto M(\tir_*(\tip,\tik),\tip,\tik)$ attains its minimum on the interval
$[0, \tipm]$. 
\end{thrm}
\begin{proof} 
Without loss of generality, assume that $\tipm \geq \sqrt{5}$.
Then for any $\tip \geq \tipm$, the function $\tir \mapsto M(\tir,\tip,\tik)$ is differentiable, so $\tir_*$ must be a root of $\dfdx{M}{\tir}(\tir,\tip,\tik)$.
We will show for $\tip$ large enough that $\tir_*$ must be a local minimum of
$\tir \mapsto H_0(\tir,\tip,\tik)$. Indeed, by Lemma \ref{lem:H0crit}, we know that if $\nu < \tip/F'_{\tip}(\sqrt{\tir_-})$, then $H_0$ has a local minimum at $\tir_* \in (0, \tir_-)$; in fact, this is also the global minimum of $\tir\mapsto M(\tir,\tip,\tik)$, since $\tir\mapsto \hat{H}(\tir,\tip,\tik)$ must be monotonically increasing on $[\tir_-,\tir_+]$ (it is increasing at both end points
and cannot have more than one critical point inside this interval). 
Thus, if we can show that $\tip/F'_{\tip}(\sqrt{\tir_-}) \to \infty$ as $\tip\to\infty$,
then we can find $\tipm = \tipm(\tik)$ such that $\nu < \tip/F'_{\tip}(\sqrt{\tir_-})$
whenever $\tip > \tipm$. 

To this end, let us show that $\tip \mapsto F'_{\tip}(\sqrt{\tir_-})$ is a positive, strictly decreasing function of $\tip$. Recall from \eqref{eq:rpm_pd} that
\[\tir_- = \tip(3+\tip - \sqrt{9+4\tip}) = \frac{2\tip^2 + \tip^3}{3+\tip + \sqrt{9+4\tip}}.\]
Substituting this into the definition of $F'_{\tip}$ and simplifying leads to
\[ F'_{\tip}(\sqrt{\tir_-}) = \frac{4\tip^2\sqrt{\tir_-}}{(\tip^2-\tir_-)^2}
= \left(\frac{(3+\sqrt{9+4\tip})^3(2+\tip)}{4\tip^2(1+\sqrt{9+4\tip})}\right)^{1/2}.\]
The latter function is positive and strictly decreasing with $\tip$; in particular, we
have $1 <  F'_{\tip}(\sqrt{\tir_-}) \leq \sqrt{8}$ whenever $\tip \geq 4$. Thus, by choosing $\tipm = \max(4, \sqrt{8}\tik)$, we guarantee that $\nu < \tip/F'_{\tip}(\sqrt{\tir_-})$ for all $\tip > \tipm$, meaning that $\tir_*(\tip,\tik)$ will always be a local minimum of $H_0$ for these values of $\tip$.

It remains to show that $H_0(\tir_*(\tip,\tik),\tip,\tik) > H_0(\tir_*(\tipm,\tik),\tipm,\tik)$ for any $\tip > \tipm$. To do so,  we write
$K(\tip,\tik) := H_0(\tir_*(\tip,\tik),\tip,\tik)$ and show that $\dfdx{K}{\tip} > 0$ for all $\tip > \tipm$.
Indeed, we have
\[
 \dfdx{K}{\tip}(\tip,\tik) = \underbrace{\dfdx{H_0}{\tir}(\tir_*(\tip,\tik),\tip,\tik)}_{=0}\dfdx{\tir_*}{\tip}(\tip,\tik) + \dfdx{H_0}{\tip}(\tir_*(\tip,\tik),\tip,\tik)
= \frac{4\tik\sqrt{\tir_*}}{\tip^2-\tir_*} > 0,
\]
since $\tir_*$ is critical point of $\tir \mapsto H_0({\tir},\tip,\tik)$ and $\tir_* < \tir_- < \tip^2$. 
The mean-value theorem then implies that $K(\tip,\tik) > K(\tipm,\tik)$ for any
$\tip > \tipm$, as claimed.

Finally, for a given $\nu > 0$, the function $\tip\mapsto M(\tir_*(\tip,\tik),\tip,\tik)$ is continuous everywhere. Since $M(\tir_*(\tip,\tik),\tip,\tik) = K(\tip,\nu)$ for all $\tip \geq \tipm$, the above result implies that
\[\inf_{\tip \geq 0}M(\tir_*(\tip,\tik),\tip,\tik) = \inf_{\tip \in [0,\tipm]} 
M(\tir_*(\tip,\tik),\tip,\tik) = \min_{\tip \in [0,\tipm]} 
M(\tir_*(\tip,\tik),\tip,\tik).\]
Thus, the minimum is attained at some optimal parameter $\tip^* \in [0,\tipm]$,
as required.
\end{proof}

The above theorem implies that $(\tir,\tip)\mapsto M(\tir,\tip,\tik)$ has a global minimum. Since $M$ is a smooth function everywhere except along the equioscillation curve $H_0(\tir,\tip,\tik) = \hat{H}(\tir,\tip,\tik)$
for ${\tip} < \sqrt{5}$, this global minimum must occur at one of the following points:
\begin{itemize}
\item A critical point of $H_0$, i.e., where $\dfdx{H_0}{\tir} = \dfdx{H_0}{\tip} = 0$;
\item A critical point of $\hat{H}$, i.e., where $\dfdx{\hat{H}}{\tir} = \dfdx{\hat{H}}{\tip} = 0$;
\item A constrainted minimum along the equioscillation curve $\tir = \rx(\tip)$, where $H_0(\tir,\tip,\tik) = \hat{H}(\tir,\tip,\tik)$. (Note that this curve is independent of $\tik$.)\\
\end{itemize}
As a result, we will need expressions for the partial derivatives of $H_0$ and $\hat{H}$, where $\dfdx{H_0}{\tir}$ is in \eqref{eq:dH0dr}, $\dfdx{\hat{H}}{\tir}$ is in \eqref{eq:dHHdr}, and the remaining two are
\begin{align}\label{eq:dH0dp}
 \dfdx{H_0}{\tip}(\tir,\tip,\tik) &= \dfdx{h}{\tip}(\sqrt{\tir},\tir,\tip,\tik) = \frac{4\tik\sqrt{\tir}}{(\tip^2-\tir)^2}(\tip^2-\tir) = \frac{4\tik\sqrt{\tir}}{\tip^2-\tir},\\
\dfdx{\hat{H}}{\tip}(\tir,\tip,\tik) &=  \dfdx{h}{\tip}(\tix_+(\tir,\tip),\tir,\tip,\tik) +\underbrace{\dfdx{h}{\tix}(\tix_+(\tir,\tip),\tir,\tip,\tik)}_{=0}\dfdx{\tix_+}{\tip}(\tix_+(\tir,\tip),\tir,\tip,\tik) \label{eq:dHHdp}\\ 
&=  \frac{4\tik\tix_+}{G(\tix_+,\tir,\tip)}(\tip^2+\tir-2\tix_+^2). \nonumber
\end{align}
Since $\dfdx{H_0}{\tip} \neq 0$ unless $\tir=0$ and $\dfdx{H_0}{\tir}\to -\infty$ as $\tir\to 0+$, we conclude that $H_0$ has no critical points at all.
We therefore focus on calculating the critical points of $\hat{H}$ and the constrained minima along the curve
$\tir = \rx(\tip)$. This will be done in the next two subsections; the results will culminate in the theorem below, which is stated in terms of the unscaled quantities
$\delta$, $p$, $T$ and $n$ and is the main theorem of this paper.
\begin{thrm}\label{thm:main_thm}
Consider the alternating OSWR method with Robin transmission conditions $\B_i = (-1)^{i+1}\partial_y + p$, where $p > 0$. Let $\tik = n\delta^2/T$,
where $\delta > 0$ is the overlap, $T > 0$ is the length of the time horizon, and $n \geq 1$ is the number of iterations taken.
\begin{itemize}
\item For large $\tik$ (e.g., for $T\to 0$, or for fixed $T$ and $\delta$ and letting $n\to\infty$), the optimized parameter $p^*$ and the corresponding $n$-step contraction factor $\rho^*$ behave asymptotically as 
$$ p^* \sim p_\infty^* := \frac{n\delta}{T} + \frac{3}{2\delta}, \qquad \rho^* \sim e^{-n^2\delta^2/T}\left(\frac{C_1n\delta^2}{T}\right)^{-n},$$
where $C_1 = 4e\approx 10.87$. 
\item For $\tik\ll 1$ (e.g., for fixed $n$ and letting either $\delta\to 0$ or $T\to\infty$), the optimized parameter $p^*$ and the corresponding $n$-step contraction factor $\rho^*$
behave asymptotically as
$$ p^* \sim p_0^* := \left(\frac{4}{9T^2\delta}\right)^{1/5}n^{2/5} + C_2\left(\frac{\delta^3}{T^4}\right)^{1/5}n^{4/5},\qquad \rho^*\sim e^{-C_3n^{6/5}(\delta^2/T)^{1/5}},$$
where $C_2 \approx 0.6346$ and $C_3 = \frac{10}{3}(\frac23)^{1/5} \approx 3.07$.
\end{itemize}
\end{thrm}

\subsection{Critical points of $\hat{H}$}
The zeros of $\dfdx{\hat{H}}{\tip}$ are characterized by
\begin{equation}\label{eq:p_stationary}
\tip^2 + \tir - 2\tix_+^2 = 0.
\end{equation}
From the definition of $\tix_+$ in \eqref{def:xstar}, we get
$ 2\tix_+^2-\tir = {\tip} + \sqrt{A\tir - B}, $
which we insert into \eqref{eq:p_stationary} to deduce
$$ \tip^2 =  {\tip} + \sqrt{A\tir - B} \iff A\tir - B = \tip^2\left(\tip-1\right)^2. $$
Using the definitions of $A$ and $B$ in \eqref{def:AB}, we can simplify the above equation to 
\begin{equation}\label{eq:px_stationary}
\tip^4 = \tir\left(\tip^2 + {2\tip}\right).
\end{equation}
We now recall that we also have $\dfdx{\hat{H}}{\tir}(\tir,\tip,\tik) = 0$, which by \eqref{eq:r_stationary} implies that
$$ 4\tik\tip\tix_+ = -{2\tip}(\tip^2 - 2\tix_+^2 - \tir) = {4\tip\tir}, $$
since $\tip^2 - 2\tix_+^2+\tir = 0$ at the critical point. Thus, at the critical point, we have $\tix_+ = \tir/\tik$, which we can again insert into \eqref{eq:p_stationary} to get
\begin{equation}
\label{eq:pr_stationary}
\tip^2 = \frac{2\tir^2}{\tik^2} - \tir = \tir\left(\frac{2\tir}{\tik^2} - 1\right).
\end{equation}
We now seek to eliminate $\tir$ from \eqref{eq:px_stationary} and \eqref{eq:pr_stationary}. First, dividing \eqref{eq:pr_stationary} by \eqref{eq:px_stationary} gives
$$\tip^{-2} = \frac{\frac{2\tir}{\tik^2}-1}{\tip^2+{2\tip}} \implies  \tir = \tik^2\left(1 + \frac{1}{{\tip}}\right).$$
Substituting into \eqref{eq:px_stationary} yields
\begin{equation}
\label{eq:prx_stationary}
\tip^4 = \tik^2\left(1 + \frac{1}{{\tip}}\right)\left(\tip^2 + {2\tip}\right) = \tik^2({\tip}+1)({\tip} + 2).
\end{equation}
We consider the case where $T$ is very small, or where we wish to take many iterations of SWR, i.e., as $\tik\to\infty$.  Letting $u=\tip/\tik$  in \eqref{eq:prx_stationary} and dividing both sides  by $\tik^4$
leads to
$$\tik^4u^4 = \tik^2(\tik u + 1)(\tik u+2) \implies u^4 = \left(u + \tik^{-1}\right)\left(u + 2\tik^{-1}\right).$$
We therefore seek a positive solution to $u^4 = (u+\epsilon)(u+2\epsilon)$ with $\epsilon := \tik^{-1}$, which approaches zero as $\tik\to\infty$.
 Since the right-hand side is greater than $u^2$, we must have
$$u^4 > u^2 \implies u^2(u^2-1) > 0 \implies u > 1. $$
Moreover, $u \to 1+$ as $\epsilon \to 0+$. In the neighbourhood of $u=1$, the implicit function theorem implies that $u$ is a differentiable function of $\epsilon$ and can be written as
$$ u = 1 + \frac32\epsilon + O(\epsilon^2). $$
Recalling that $\tip = \tik u$, we arrive at an asymptotic expansion of $\tip$ for large $\tik$, which in turn yields those for $\tir$ and $\tix$ at this critical point $(\tir^*,\tip^*)$, namely:
\begin{align}
 \tip^* &= \tik \left(1 + \frac{3}{2}\epsilon + O(\epsilon^2)\right), \label{pinf_asympt} \\
\tir^* &= \tik^2\left(1+\epsilon + O(\epsilon^{2})\right),\\
\tix_+ &= \tik\left(1+\epsilon + O(\epsilon^{2})\right).
\end{align}
We are now ready to calculate the asymptotic expression for $\hat{H}(\tir,\tip,\tik)$ at the critical point $(\tir^*,\tip^*)$. Since
$$ \hat{H}(\tir^*,\tip^*,\tik) = h(\tix_+,\tir^*,\tip^*,\tik) = \tir^* - 2\tik \tix_+ + \tik\log\left[\frac{(\tip^*-\tix_+)^2+\tix_+^2-\tir^*}{(\tip^*+\tix_+)^2+\tix_+^2-\tir^*}\right],$$
we compute the relevant parts of the expression as follows:
\begin{align*}
\tir^* - 2\tik \tix_+ &= -\tir^* = -\tik^2(1+\epsilon + O(\epsilon^2)),\\
\tip^*-\tix_+ &= \tik\epsilon \left(\frac12 + O(\epsilon) \right) = \frac12+O(\epsilon),\\
\tip^*+\tix_+ &= \tik \left(2 + \frac52\epsilon + O(\epsilon^2)\right),\\
\tix_+^2-\tir &= \tix_+(\tix_+-\tik) = \tik^2\epsilon(1+\epsilon+O(\epsilon^2)) = \tik(1+\epsilon+O(\epsilon^2)),\\
\frac{(\tip^*-\tix_+)^2+\tix_+^2-\tir^*}{(\tip^*+\tix_+)^2+\tix_+^2-\tir^*} &= 
\frac{\tik(1+\frac54\epsilon + O(\epsilon^2))}{\tik^2(4 + 11\epsilon + O(\epsilon^2))} = \frac14\epsilon \left(1-\frac32\epsilon  + O(\epsilon^2)\right),\\
\hat{H}(\tir^*,\tip^*,\tik) &= -\tik^2(1+\epsilon + O(\epsilon^2)) + \tik \log(\epsilon/4) - \tik\left(\frac32\epsilon + O(\epsilon^2)\right)\\ &= -\tik^2\left(1+\epsilon|\log\epsilon| + \epsilon(1+\ln 4) + O(\epsilon^2)\right).
\end{align*}
In order to see whether $\hat{H}(\tir^*,\tip^*,\tik)$ yields a correct convergence bound, we need to compare this value with $H_0(\tir^*,\tip^*,\tik)$, which we can calculate by
simply replacing $\tix_+$ by $\sqrt{\tir^*}$ in the expression for $\hat{H}$:
\begin{align*}
\tir^* - 2\tik\sqrt{r^*} &= -\tik^2(1+O(\epsilon^2)),\\
\tip^*-\sqrt{\tir^*} &= \tik\epsilon(1+O(\epsilon)),\\
\tip^*+\sqrt{\tir^*}&= \tik(2+2\epsilon+O(\epsilon^2)),\\
\frac{\tip^*-\sqrt{\tir^*}}{\tip^*+\sqrt{\tir^*}} &= \frac{\epsilon}{2} + O(\epsilon^2)\\
H_0(\tir^*,\tip^*,\tik) &= -\tik^2(1+O(\epsilon^2)) - 2\tik|\log(\epsilon/2 + O(\epsilon^2))| \\&= -\tik^2(1+2\epsilon|\log\epsilon| + O(\epsilon)).
\end{align*}
Since $H_0(\tir^*,\tip^*,\tik) < \hat{H}(\tir^*,\tip^*,\tik)$ asymptotically, we conclude that $\hat{H}(\tir^*,\tip^*,\tik)$ is a valid convergence bound. Note that by choosing this $p^* = \delta \tip^*$
as the Robin parameter, we obtain an $n$-step contraction factor of the form
$$\rho(r^*,p^*,n) = \exp\left(-\frac{T}{\delta^2}(\tik^2 + \tik\log \tik + \tik(1+\ln 4) + O(1))\right) \simeq e^{-n^2\delta^2/T}\left(\frac{Cn\delta^2}{T}\right)^{-n},$$
where $C = 4e \approx 10.87$. We therefore gain an extra factor of $\left(\frac{Cn\delta^2}{T}\right)^{-n}$ relative to the Dirichlet bound,
which indicates faster convergence than for Dirichlet transmission conditions.\\

A different scaling regime consists of letting $\tik\to 0$, i.e., fixing the number of iterations $n$ and considering what happens in the limit of a vanishing overlap $\delta \to 0$ (or, equivalently, that of a long time horizon $T\to \infty$). In this case, making the substitution $\tip = \tik^{1/2}u$ in \eqref{eq:prx_stationary} leads to
$$ u^4 = (\tik^{1/2}u+1)(\tik^{1/2}u+2) = 2 + \epsilon(3u + \epsilon u^2), \qquad \text{where $\epsilon = \tik^{1/2}$.} $$
Taking the limit as $\tik\to 0$ leads to $u^4 = 2$, whose only positive solution is $\bar{u} = 2^{1/4}$. 
By Rouch\'e's theorem, \eqref{eq:prx_stationary} will possess a simple root near
$\bar{u}$ for $\tik$ small enough; moreover, the trajectory of this root is differentiable with respect to $\epsilon$, since $\bar{u}$ is a simple root of the limiting equation.
Thus, to obtain the asymptotics, we write $u^*= \bar{u} + c_1\epsilon + O(\epsilon^2)$ and insert it into \eqref{eq:prx_stationary}. Solving for $c_1$ leads to
$$ u^* = 2^{1/4} + \frac{3}{2^{5/2}}\epsilon + O(\epsilon^2). $$
We can now compute the asymptotic expansions of $\tip^*$, $\tir^*$ as before. This leads to
$$
\tip^* = \tik^{1/2}\left(2^{1/4} + \frac{3}{2^{5/2}}\tik^{1/2} + O(\tik)\right),\qquad
\tir^* = \tik^{3/2}\left(2^{-1/4} - \frac58\tik^{1/2} + O(\tik)\right).
$$
We now continue the calculation to obtain the asymptotics for $\hat{H}(\tix_+,\tir^*,\tip^*,\tik)$. Note that the coefficients are given as floating
point numbers in order to facilitate comparisons with ${H}_0(\tix_+,\tir^*,\tip^*,\tik)$ later:
\begin{align*}
\tir^* - 2\tik\tix_+ &= -\tir^* = -\tik^{3/2}\left(0.8409 - 0.625\tik^{1/2} + O(\tik)\right),\\
\tip^*-\tix_+ &= \tik^{1/2}\left(0.3483 +  1.1553\tik^{1/2} + O(\tik) \right),\\
\tip^*+\tix_+ &= \tik^{1/2} \left(2.0301 - 0.0947\tik^{1/2} + O(\tik)\right),\\
\tix_+^2-\tir &= \tik\left(0.707-1.892\tik^{1/2} + O(\tik)\right),\\
\frac{(\tip^*-\tix_+)^2+\tix_+^2-\tir^*}{(\tip^*+\tix_+)^2+\tix_+^2-\tir^*} &= 
\frac{\tik(0.8283 - 1.0872\tik^{1/2} + O(\tik))}{\tik(4.8283 - 2.2765\tik^{1/2} + O(\tik))} = 0.1716 + 0.2463\tik^{1/2} + O(\tik),\\
\hat{H}(\tir^*,\tip^*,\tik) &= -\tik^{3/2}(0.8409 + O(\tik^{1/2})) + \tik\left(\log(0.1716) +  1.4353\tik^{1/2} + O(\tik)\right)\\
&= -1.7626\tik + 0.5944\tik^{3/2} + O(\tik^2).
\end{align*}
As before, we compare the above value with $H_0(r^*,p^*,k)$:
\begin{align*}
\tir^* - 2\tik\sqrt{\tir^*} &= \tik^{3/2}\left(2^{-1/4} + O(\tik^{1/4})\right),\\
\tip^*\pm \sqrt{\tir^*} &= \tik^{1/2}2^{1/4}\left(1\pm 2^{-3/8}\tik^{1/4} + 3\cdot 2^{-11/4}\tik^{1/2} + O(\tik^{3/4})\right),\\
H_0(\tir^*,\tip^*,\tik) &= -3.0844\tik^{5/4} + O(\tik^{7/4}).
\end{align*}
This time we have $H_0(\tir^*,\tip^*,\tik) > \hat{H}(\tir^*,\tip^*,\tik)$ whenever $\tik = n\delta^2/T \lesssim 0.0728$, meaning that this critical point does not yield a valid convergence bound as
$\tik \to 0$. We need to look instead for a constrained minimum along the equioscillation curve to obtain a valid bound in this case.
\subsection{Constrained minimum along the equioscillation curve}
We start by noting that the equioscillation curve $\rx = \rx(\tip)$ is \emph{independent of $\tik$}, since 
\begin{equation}
\label{def:gamma}
\begin{aligned}
\gamma(\tir,\tip) &:=
\tik^{-1}(\hat{H}(\tir,\tip,\tik)-H_0(\tir,\tip,\tik)) \\
&= 2(\sqrt{\tir}-\tix_+) + \log\left(\frac{(\tip-\tix_+)^2+\tix_+^2-\tir}{(\tip+\tix_+)^2+\tix_+^2-\tir}\right) - \log\left(\frac{(\tip-\sqrt{\tir})^2}{(\tip+\sqrt{\tir})^2}\right)
\end{aligned}
\end{equation}
is independent of $\tik$ (including $\tix_+ = \tix_+(\tir,\tip)$). However, the constrained minimum will depend on $\tik$, since we are minimizing the function $H_0(\tir,\tip,\tik)$ there, and the latter depends on $\tik$. Our strategy, as before, is to find the correct scaling that reveals the asymptotic behaviour of $\tip$, $\tir$ and $\tix_+$ as $\tik\to 0$. To do so, we will need the following
lemma, whose proof appears in Appendix \ref{app:pf_lem4.3}.

\begin{lmm}\label{thm:equi_curve}
Let $\tip \in (0, \sqrt{5})$ be given, and let $\rx \in [0,\tir_-(\tip)]$ be such that $0 \leq \rx \leq \tir_-$ and $H_0(\rx,\tip,\tik) =\hat{H}(\rx,\tip,\tik)$ (cf.~Theorem \ref{thm:M_prop}). Then there exist positive constants $C_1$ and $C_2$, independent of $\tip$, such that $C_1 \tip^3\leq \rx \leq C_2\tip^3$.
\end{lmm} 

Motivated by the above lemma, we make the substitution $\rx = c^2\tip^3$ and $\tix_+ = \xi\sqrt{\tip}$ into $\gamma(\tir,\tip) = 0$ to get
\begin{equation}
\label{eq:pxi1}
\begin{aligned}
 0 &= 2(c\tip^{3/2} - \xi \tip^{1/2}) + \log\left(\frac{\tip[(\tip^{1/2} - \xi)^2 + \xi^2-c^2\tip^2]}{\tip[(\tip^{1/2} + \xi)^2 + \xi^2-c^2\tip^2]}\right) - \log\left(\frac{\tip^2(1-c\tip^{1/2})^2}{\tip^2(1+c\tip^{1/2})^2}\right)\\
&= 2\tip^{1/2}(c\tip - \xi) + \log\left(\frac{(\tip^{1/2} - \xi)^2 + \xi^2-c^2\tip^2}{(\tip^{1/2} + \xi)^2 + \xi^2-c^2\tip^2}\right) - \log\left(\frac{(1-c\tip^{1/2})^2}{(1+c\tip^{1/2})^2}\right).
\end{aligned}
\end{equation}
Additionally, the definition of $\tix_+$ in \eqref{def:xstar} leads to the relation
\begin{align}\nonumber
2\tix_+^2 = 2\xi^2\tip = c^2\tip^3 + \tip + \sqrt{(4\tip+2\tip^2)c^2\tip^3 -(\tip^4+2\tip^3-\tip^2)},\\
\intertext{which implies}
2\xi^2 = 1 + c^2\tip^2 + \sqrt{c^2(4\tip^2+2\tip^3) + (1-2\tip-\tip^2)}.\label{eq:pxi2}
\end{align}
We wish to find the simultaneous solution $(c,\xi)$ of \eqref{eq:pxi1}--\eqref{eq:pxi2}, and particularly in the limit of $\tip\to 0$, since $\tip^*\to 0$ as $\tik\to 0$. Indeed,
$\tip^*$ must be bounded because we are in the regime where $0 < \tip < \sqrt{5}$ (otherwise the equioscillation curve
plays no role in the minimization problem). If the optimal $\tip$ approaches a positive constant as $\tik\to 0$, then $\rx$ and $\tix_+$ would also approach a positive constant,
meaning that $H_0(\tir, \tip, \tik)$ approaches a constant as well. This would lead to an asymptotic bound that does not contain $\tik$, i.e., it does not contain the iteration number $n$.
Such a bound would be useless in practice.

Letting $\tip\to 0$ in \eqref{eq:pxi2} leads to $\xi =1$. Next, we rewrite \eqref{eq:pxi1} as
\begin{equation}
\label{eq:pxi3}
\xi-c\tip = \frac{1}{2\tip^{1/2}}\log\left(\frac{(\tip^{1/2} - \xi)^2 + \xi^2-c^2\tip^2}{(\tip^{1/2} + \xi)^2 + \xi^2-c^2\tip^2}\right) -  \frac{1}{2\tip^{1/2}}\log\left(\frac{(1-c\tip^{1/2})^2}{(1+c\tip^{1/2})^2}\right).
\end{equation}
Taking the limit $\tip\to 0$ on both sides yields
$$\xi = -\frac{1}{\xi} + 2c \implies c = 1. $$
We therefore conclude that $\rx = \tip^3 + O(\tip^4)$. Moreover, after dividing \eqref{eq:pxi1} by $\tilde{p}^{1/2}$, we see that the Jacobian of the system \eqref{eq:pxi1}--\eqref{eq:pxi2} with respect to $(c,\xi)$ is invertible in the neighbourhood of $(c,\xi,\tip) = (1,1,0)$; thus, $c$ and $\xi$ are locally differentiable (in fact, analytic) with respect to $\tip$. Substituting the Taylor series of these functions into \eqref{eq:pxi1}--\eqref{eq:pxi2} and equating coefficients leads to
$$ \xi = 1 - \frac14 \tip + \frac{15}{32}\tip^2 + O(\tip^3),\qquad c = 1 - \frac{11}{12}\tip + \frac{219}{160}\tip^2 + O(\tip^3). $$
We are now ready to prove the asymptotic optimal parameter $\tip$ along the equioscillation curve.
\begin{thrm}
Let $\tip \in (0,\sqrt{5})$ and $\rx = \rx(\tip)$ be the equioscillation curve satisfying $H_0(\rx,\tip,\tik) = \hat{H}(\rx,\tip,\tik)$. Then the minimum of $\tip\mapsto H_0(\rx(\tip),\tip,\tik)$
has the asymptotic behaviour
$$ \tip^* =\Bigl(\frac{2}{3}\Bigr)^{2/5}\nu^{2/5} + K\nu^{4/5} + O(\nu^{6/5}),\qquad H_0(\rx(\tip^*),\tip^*,\tik)\sim -C\,\tik^{6/5}$$
with $K\approx 0.6346$ and $C =\frac{10}{3}\left(\frac23\right)^{1/5}\approx 3.07$, where we write $f \sim \alpha \tik^\beta$ whenever $f = \alpha \tik^\beta + o(\tik^{\beta})$ as $\tik\to 0$.
\end{thrm}
\begin{proof}
Let $V(\tip) = H_0(\rx(\tip),\tip,\tik)$. We seek the value of $\tip$ minimizing $V$ by finding a root of $V'(\tip)$. From \eqref{eq:dH0dr} and
\eqref{eq:dH0dp}, we see that
\begin{align*}
V'(\tip) &= \dfdx{H_0}{\tip}(\rx(\tip),\tip,\tik) + \dfdx{H_0}{\tir}(\rx(\tip),\tip,\tik)\frac{d\rx}{d\tip}\\
&= \left[\frac{4\tik\sqrt{\rx}}{\tip^2-\rx}+\left(1-\frac{\tik}{\sqrt{\rx}} - \frac{2\tik\tip}{\sqrt{\rx}(\tip^2-\rx)}\right)\frac{d\rx}{d\tip} \right].
\end{align*}
Setting $V'(\tip) = 0$ leads to
$$ \frac{\tik}{\sqrt{\rx}(\tip^2-\rx)}\left(((\tip^2-\rx) + 2\tip)\frac{d\rx}{d\tip} - 4\rx\right) = \frac{d\rx}{d\tip},$$
leading to
\begin{equation}\label{eq2}
 \tik = \frac{\sqrt{\rx}(\tip^2-\rx) (\rx)'}{((\tip^2-\rx) + 2\tip)(\rx)' - 4\rx}.
\end{equation}
Since $\rx = c^2\tip^3$, we have the asymptotic expansions
\begin{align*}
 \sqrt{\rx} &= c\tip^{3/2} = p^{3/2}\left(1 - \frac{11}{12}\tip + \frac{219}{160}\tip^2 + O(\tip^3)\right),\\
(\rx)' &= 3\tip^2 - \frac{22}{3}\tip^3 + \frac{161}{9}\tip^4 + O(\tip^5).
\end{align*}
Plugging these into \eqref{eq2} yields 
$$\nu = \tip^{5/2}\left(\frac32 - \frac{79}{24}\tip + \frac{3841}{576}\tip^2 + O(\tip^3)\right). $$
To express $\tip$ as a function of $\nu$, we first let $\tip = q^2\nu^{2/5}$, so that
$$ \nu = q^5\nu\left(\frac32 - \frac{79}{24}q^2\nu^{2/5}+ \frac{3841}{576}q^4\nu^{4/5} + O(q^6\nu^{6/5})\right).$$
Now we divide both sides by $\nu$ and let $z = \nu^{2/5}$ to get
$$ 1 = q^5\left(\frac32 - \frac{79}{24}zq^2+ \frac{3841}{576}z^2q^4 + O(z^3q^6)\right). $$
This equation is now analytic in $z$ and $q$ and has a simple root at $(q,z) = ((2/3)^{1/5},0)$. The implicit function theorem 
once again allows us to obtain an asymptotic expansion of $q$ in terms of $z = \nu^{2/5}$, namely
\begin{alignat}{2}
q &= \Bigl(\frac{2}{3}\Bigr)^{1/5} + q_1\nu^{2/5}  + O(\nu^{4/5}),&&\qquad \text{where $q_1\approx 0.3441$.}\nonumber
\intertext{This finally leads to}
\tip &=\nu^{2/5}\left(\Bigl(\frac{2}{3}\Bigr)^{2/5} + K\nu^{2/5} + O(\nu^{4/5})\right),&&\qquad \text{where $K \approx 0.6346$.}\label{p0_asympt}
\end{alignat}
To calculate the asymptotic contraction rate, we substitute $\tip^* \sim (\frac23\tik)^{2/5}$ into $H_0$ to obtain
\begin{align*}
H_0(\rx,\tip^*,\tik) &= \rx-2\tik\sqrt{\rx}+\tik\log\left(\frac{(\tip-\sqrt{\rx})^2}{(\tip+\sqrt{\rx})^2}\right)\\
&\sim \tip^3 - 2\tik\tip^{3/2} +2\tik\log\left(\frac{1-\tip^{1/2}}{1+\tip^{1/2}}\right)\\
&\sim \tip^3 - 2\tik\tip^{3/2} - 4\tik\tip^{1/2}\\
&\sim \left(\frac{2\tik}{3}\right)^{6/5} + O(\tik^{8/5}) - 4\left(\frac23\right)^{1/5}\tik^{6/5} \sim -C\tik^{6/5},
\end{align*}
where $C =4(\frac23)^{1/5} - (\frac23)^{6/5}   = \frac{10}{3}(\frac23)^{1/5}$.
\end{proof}

\begin{figure}
\begin{center}
\includegraphics[width=0.48\textwidth,bb=110 235 490 545,clip]{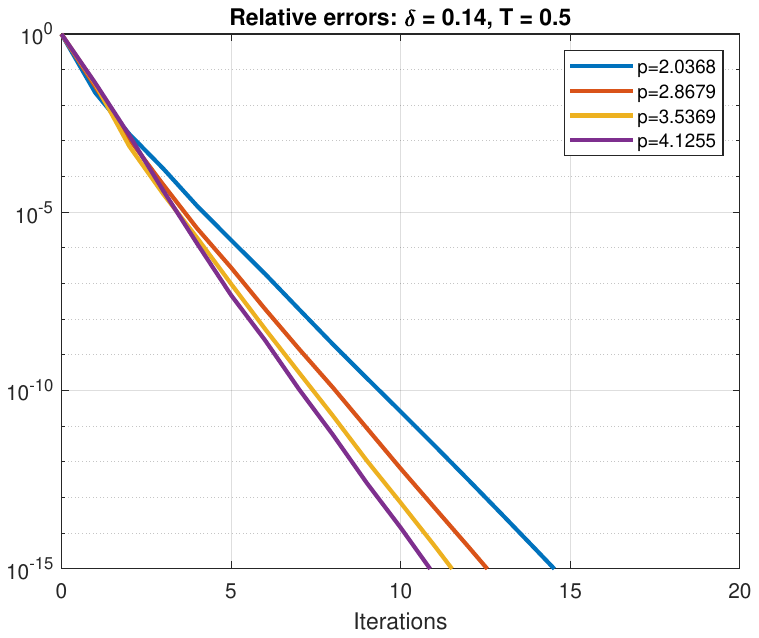}
\includegraphics[width=0.48\textwidth,bb=110 235 490 545,clip]{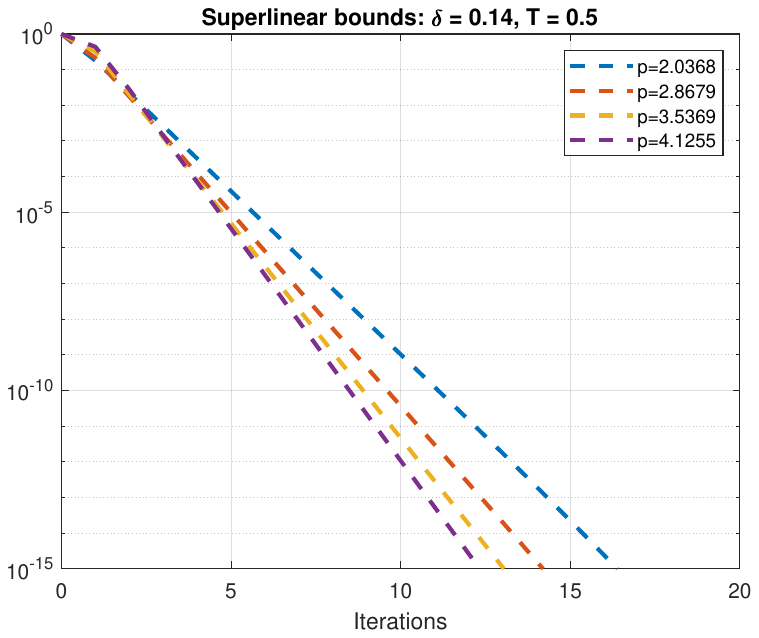}\\
\includegraphics[width=0.48\textwidth,bb=110 235 490 545,clip]{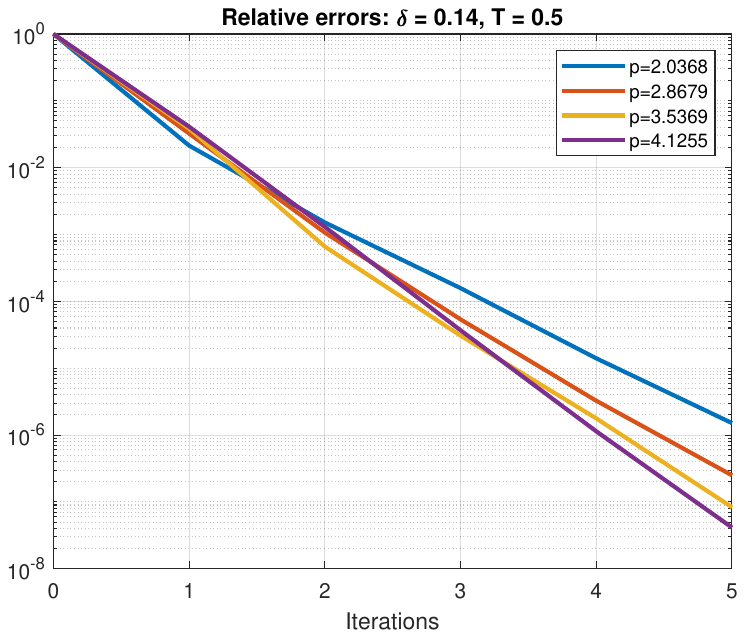}
\includegraphics[width=0.48\textwidth,bb=110 235 490 545,clip]{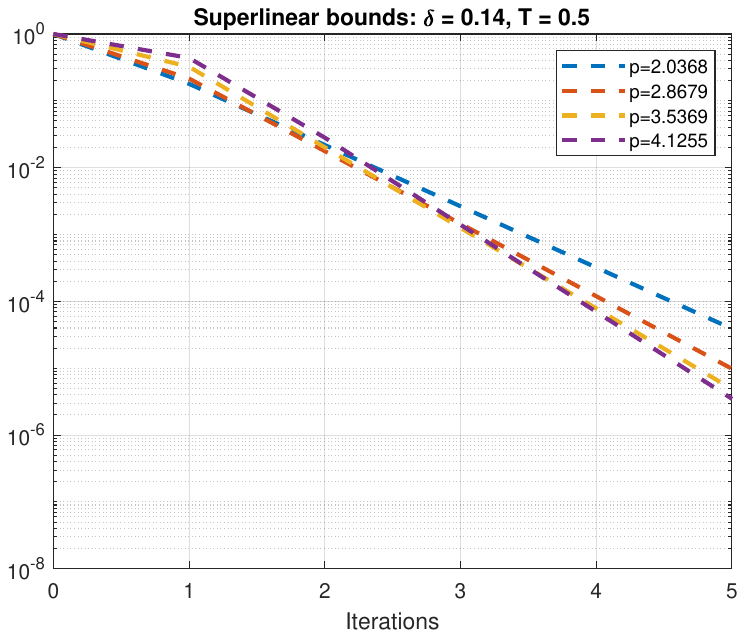}\\
\includegraphics[width=0.48\textwidth,bb=110 235 490 545,clip]{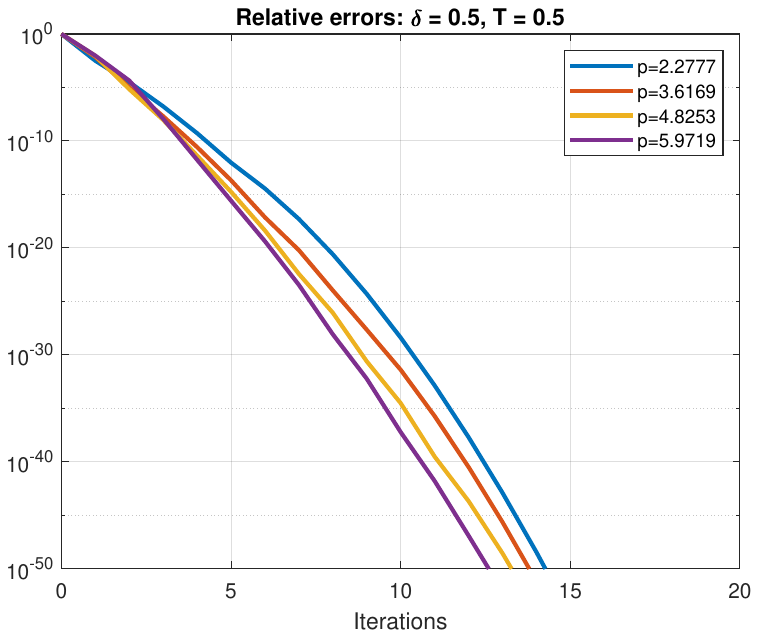}
\includegraphics[width=0.48\textwidth,bb=110 235 490 545,clip]{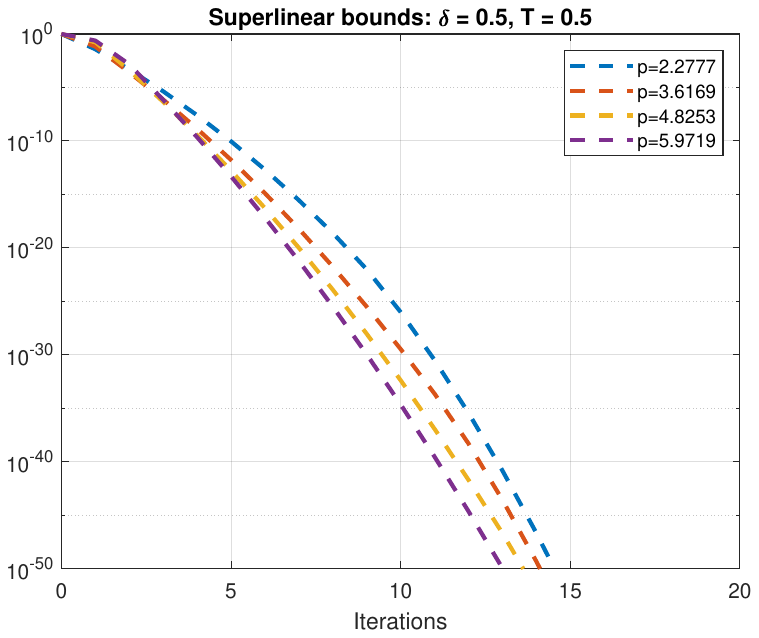}\\
\includegraphics[width=0.48\textwidth,bb=110 235 490 545,clip]{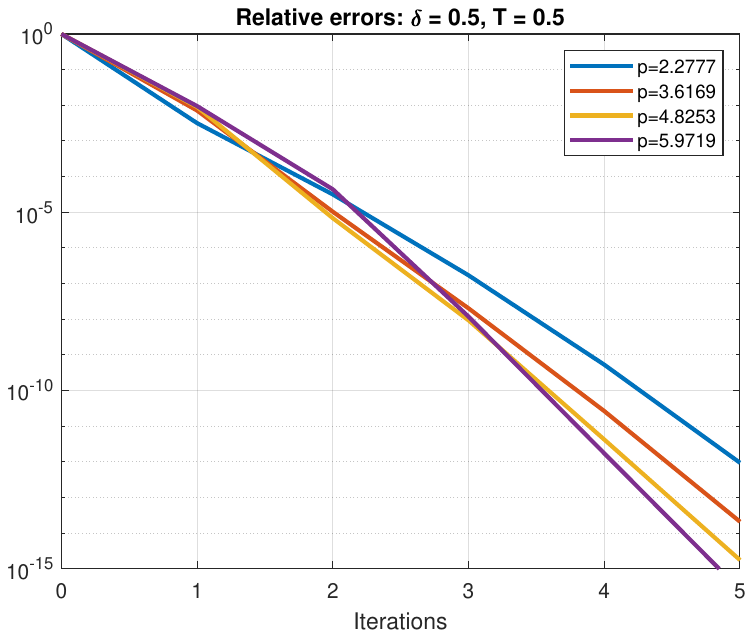}
\includegraphics[width=0.48\textwidth,bb=110 235 490 545,clip]{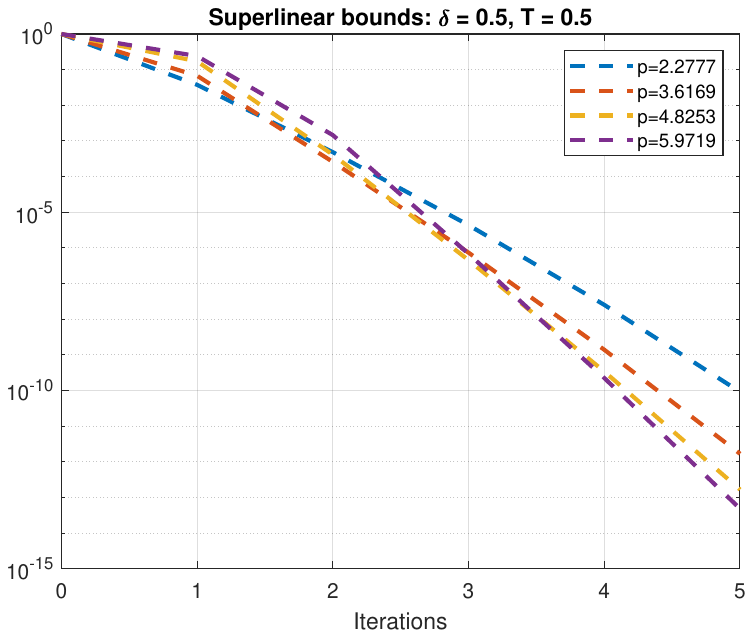}\\
\caption{Relative $L^2$ errors of OSWR iterations and their theoretical bounds for $T=0.5$ and $\delta\in \{0.14, 0.5\}$, with the optimal parameter $p$ chosen for $n=1,2,3,4$ in each case. Left: relative errors, right: superlinear bounds. The second and fourth rows are zooms of the first and third ones to better show
the behaviour of the first 5 iterations.\label{fig:rel_errs}}
\end{center}
\end{figure}
\section{Numerical experiments}
We start by studying the convergence of OSWR using optimized parameters $p^* = p^*(n)$ that minimize the bound $M(r^*(p),p,n)$ for $n=1,2,3,4$; we consider $T=0.5$ and two
overlap sizes, $\delta \in \{0.14, 0.5\}$. In each case, we show two sets of curves: the $L^2$ relative error of OSWR iterations applied to the same discretized problem
as in Section \ref{sec:3.3}, i.e., the heat equation $u_t=u_{yy}$ with Dirichlet boundary conditions on the finite domain $\Omega = (-4,4)$, with centred difference in space and backward Euler in time with $\Delta y = 1/800$
and $\Delta t = 1/2400$. We then plot, next to the convergence curves, the theoretical bounds for the same $p^*$, see Figure \ref{fig:rel_errs}. 
As we mentioned before, the optimal $p^*$ depends on $n$, the target number of iterations, via the parameter $\tik=n\delta^2/T$; 
indeed, we see from Figure \ref{fig:rel_errs} that while we have
convergence for any choice of parameter $p$, the best parameter does depend on the target number of iterations, as seen by the curves intersecting each other
during the first four iterations of OSWR.

\begin{rmrk}
For discretized problems, convergence is technically always geometric, since the problem is linear and finite dimensional;
we can only have superlinear convergence if the iteration is nilpotent. Otherwise, the spectral radius is always strictly positive and corresponds to
the smallest eigenvalue of the iteration matrix, which shows up as a pole in Laplace space. Thus, the superlinear bounds we obtain from the continuous setting are technically invalid for discretized problems. Nonetheless, when the overlap between the subdomains are wide enough (i.e., when it spans at least a few mesh widths),
the spectral radius is so small that the superlinear bounds remain valid in practice. This is why we use a very fine discretization in all our experiments, so that
 $\delta \gg \Delta y$ in all cases. The behaviour of OSWR when $\delta \approx \Delta y$ will be studied in a subsequent paper.
\end{rmrk}

To illustrate the asymptotic behaviour of $p^*$ stated in Theorem \ref{thm:main_thm}, we report in Table \ref{table1} the computed $p^*$ as well as the asymptotic values $p^*_\infty$ and $p_0^*$ for the same value of $\tik$. We also report whether the minimum is characterized by equioscillation, rather than by a critical point of $\hat{H}$. We observe for $\delta = 0.14$ that $\frac{|p^* - p^*_0|}{p^*} \approx C\tik^{4/5}$  since $\tik$ is small;
this is consistent with the asymptotic expansion \eqref{p0_asympt}. For $\delta = 0.5$, we see that we are in a transition zone, with $p^*$ closer to $p_0^*$ for $n=1$,
but closer to $p_\infty^*$ for $n=4$. Although \eqref{pinf_asympt} suggests that $p_\infty^*$ should be an $O(\nu^{-2})$ approximation of $p^*$, 
our $\nu$ values are too far from the asymptotic regime of $\nu\to\infty$, so we only observe an $O(\nu^{-1})$ approximation in this case. To illustrate \eqref{pinf_asympt} properly, we show results for $\delta=1$, 
which yields larger values of $\nu$; here, we see that $p_\infty^*$ is a much better approximation of $p^*$, and the relative approximation error now scales like $O(\nu^{-2})$,
as expected. 

\begin{figure}
\includegraphics[width=0.46\textwidth,bb=0 -25 702 560,clip]{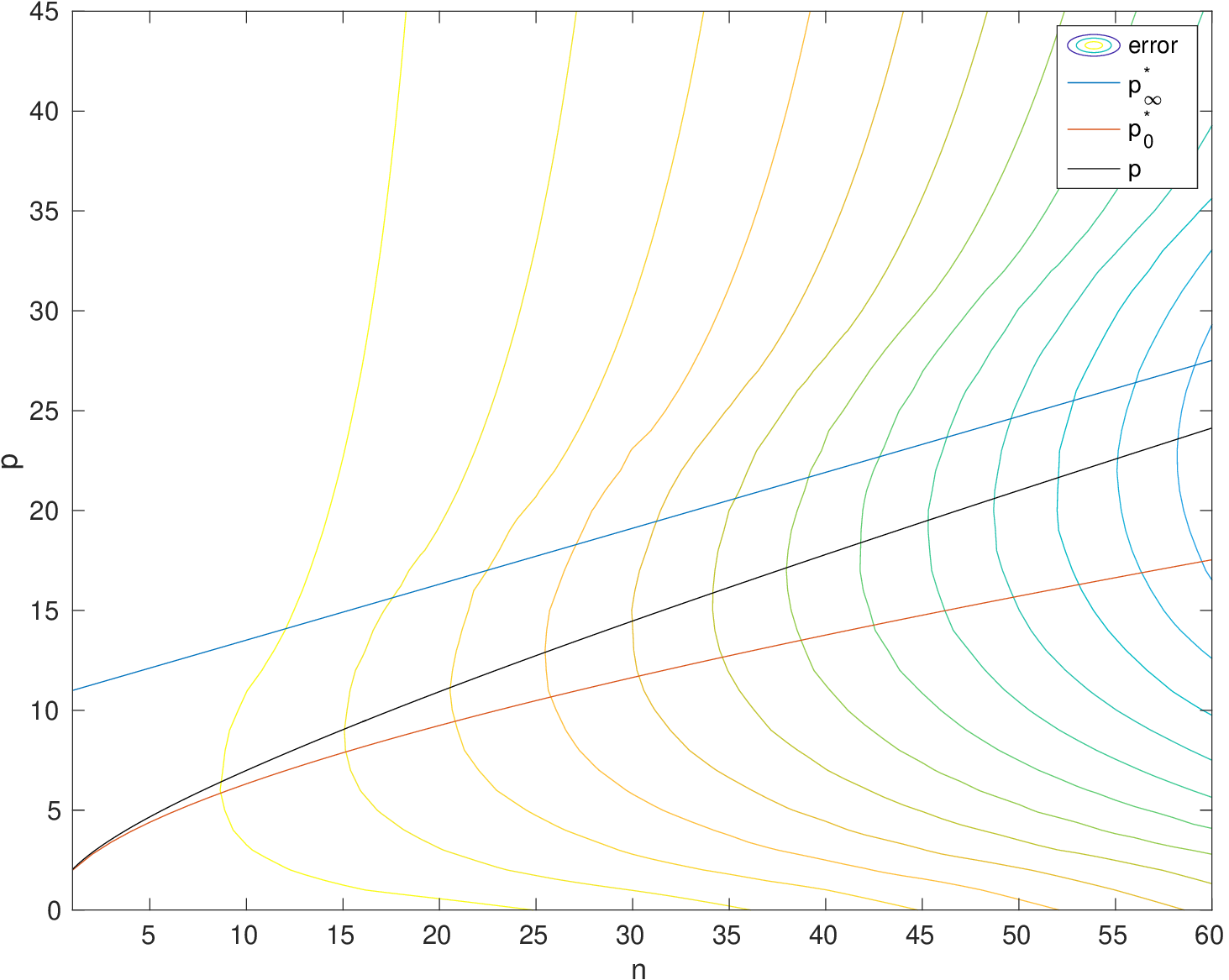}
\includegraphics[width=0.49\textwidth,bb=110 240 485 537,clip]{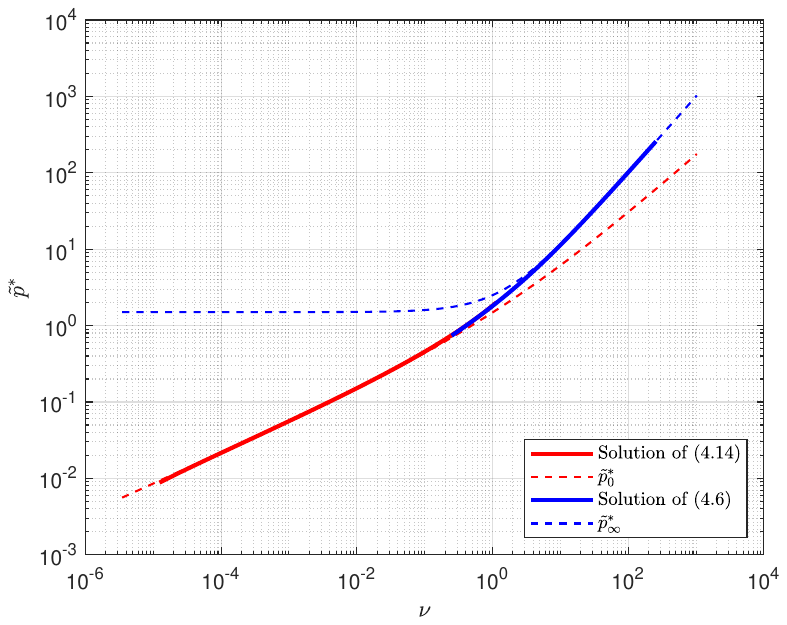}
\caption{Left: The best parameter $p$ (black curve) calculated by solving \eqref{eq:prx_stationary} or \eqref{eq2} for $T=0.5$, $\delta=0.14$ and various $n$.
The asymptotic formulas $p_0^*$ and $p_\infty^*$ are shown in red and blue, and the iteration errors are shown as contours in the background. Right: comparison
between the best parameter $p^*$ (solid) and the asymptotic formulas (dashed) for various $\tik$.\label{fig6}}
\end{figure}

In the left panel of Figure~\ref{fig6}, we show the level sets of the error
measured in the space-time infinity norm when applying OSWR to the
error equations for the configuration considered in the top part of
Table~\ref{table1}. We vary the number of iterations $n$ and the Robin parameter
$p$, using the same random initial guess in each case. The error
decreases with increasing $n$ for all values of $p$. The optimal
parameter $p$, computed by solving \eqref{eq:prx_stationary} or \eqref{eq2}, is shown in black
and closely follows the valley of minimal error indicated by the level
sets. We also plot the asymptotic approximations $p_0^*$ and
$p_{\infty}^*$ in red and blue. These approximations capture the
behaviour of the optimal parameter well: $p_0^*$ follows the valley for
small $n$, while $p_{\infty}^*$ will provide a better approximation
when $n$ becomes large.

To illustrate the behaviour of the black curve (the valley of minimal error) over a wider range of $\nu$, we plot
in the right panel of Figure~\ref{fig6} the optimal parameter $\tilde{p}^*$ as a function of $\nu$; we further distinguish
the portion of the curve given by \eqref{eq:prx_stationary} 
(critical point of $\hat{H}$, blue solid curve) from the part given by \eqref{eq2} (equioscillation, red solid curve).
We see that the two parts of the curve connect very smoothly in the transition zone, so solving either equation will yield a good approximation for $p^*$. Moreover, we show in 
dashed lines the asymptotic formulas $p_0^*$ and $p_\infty^*$. 
We see that $p_0^*$ is a good approximation of $p^*$ for $\nu\lesssim 0.5$, whereas $p_\infty^*$ approximates $p^*$ well for $\nu \gtrsim 5$. Since neither
approximation is accurate enough for $\nu \in (0.5,5)$, we recommend solving \eqref{eq:prx_stationary}, namely the quartic equation $\tip^4 = \nu^2(\tip+1)(\tip+2)$, to obtain
$p^*$ directly for such values of $\nu$.

\begin{table}[t]
\centering
\caption{Optimal values of $p$ for different numbers of iterations, compared with their asymptotic values $p^*_\infty$ and $p^*_0$ given by Theorem \ref{thm:main_thm}. 
We also indicate whether the optimum lies on the equioscillation curve and calculate the discrepancy with one of the asymptotic values. The time horizon is $T=0.5$ for all cases. \label{table1}}
\begin{tabular}{|c|c|c|c|c|c|c|}
\hline
\multicolumn{7}{|c|}{$\delta = 0.14$}\\
\hline
$n$ & $\tik$ & $p^*$ & $p_\infty^*$ & $p_0^*$ & Equiosc.? & $\frac{|p^*-p_0^*|}{p^*\tik^{4/5}}$\\
\hline
1 & 0.0392 & 2.0368 & 10.9943 & 2.0021 & Yes & 0.2275\\
2 & 0.0784 & 2.8679 & 11.2743 & 2.7849 & Yes & 0.2217\\
3 & 0.1176 & 3.5369 & 11.5543 & 3.3978 & Yes & 0.2180\\
4 & 0.1568 & 4.1255 & 11.8343 & 3.9241 & Yes & 0.2150\\
\hline
\multicolumn{7}{|c|}{$\delta = 0.5$}\\
\hline
$n$ & $\tik$ & $p^*$ & $p_\infty^*$ & $p_0^*$ & Equiosc.? & $\frac{|p^*-p_\infty^*|}{p^*\tik^{-1}}$\\
\hline
1 & 0.5 & 2.2777 & 4.0 & 2.0178 & Yes & 0.3781\\
2 & 1.0 & 3.6169 & 5.0 & 2.9698 & No  & 0.3824\\
3 & 1.5 & 4.8253 & 6.0 & 3.7555 & No  & 0.3652\\
4 & 2.0 & 5.9719 & 7.0 & 4.4537 & No  & 0.3443\\
\hline
\multicolumn{7}{|c|}{$\delta = 1$}\\
\hline
$n$ & $\tik$ & $p^*$ & $p_\infty^*$ & $p_0^*$ & Equiosc.? & $\frac{|p^*-p_\infty^*|}{p^*\tik^{-2}}$\\
\hline
1 & 2.0 & 2.9860 & 3.5 & 2.2269 & No & 0.6885\\
2 & 4.0 & 5.1504 & 5.5 & 3.4042 & No  & 1.0861\\
3 & 6.0 & 7.2325 & 7.5 & 4.4020 & No  & 1.3315\\
4 & 8.0 & 9.2827 & 9.5 & 5.3029 & No  & 1.4982\\
\hline
\end{tabular}
\end{table}

\section{Conclusion}
We have developed a new analytical technique for capturing the
superlinear convergence behaviour of Schwarz waveform relaxation
methods. Traditionally, analyzing this regime has relied on
technically demanding kernel estimates, which have hindered the
optimization of transmission conditions in optimized variants of these
algorithms for almost three decades. Yet, the superlinear convergence
regime is precisely where these methods achieve their best
performance, making the optimization of transmission conditions
particularly important.

Our approach exploits the real part of the Laplace parameter in a
carefully chosen manner to obtain sharp convergence estimates without
resorting to kernel techniques. This enables us, for the first time
since the introduction of this class of algorithms three decades ago,
to characterize in detail how the optimal transmission parameter
depends on the number of iterations, the overlap size, and the time
window length. Furthermore, our analysis yields simple asymptotic
formulas that can be readily used in practice to select effective
transmission parameters.

We have focused on the one-dimensional heat equation to develop and
validate our approach. However, the underlying analytical technique is
not restricted to this setting and can be extended to more general
problems, which we leave for future work.

\appendix
\section{Proof of Lemma \ref{thm:equi_curve}}\label{app:pf_lem4.3}
First, we need a bound for $\tir_-$. By definition, we recall from \eqref{eq:rpm_pd} that 
$$ \tir_- = {\tip}\left(3+{\tip} - \sqrt{9 + 4{\tip}}\right)
= {\tip}\left(\frac{2{\tip} + \tip^2}{3+{\tip} + \sqrt{9 + 4{\tip}}}\right) \leq \tip^2\left(\frac{2+{\tip}}{3+{\tip}+\sqrt{9+4{\tip}}}\right).$$
Since the term inside the parentheses is an increasing function of ${\tip}$ (verified by differentiation), the fact that ${\tip} < \sqrt{5}$ implies
\begin{equation}
\label{eq:r_is_small}
\frac{\tir_-}{\tip^2} <\frac{2+\sqrt{5}}{3+\sqrt{5}+\sqrt{9+4\sqrt{5}}}  = \frac{1}{\sqrt{5}}.
\end{equation}
Since $\rx \leq \tir_-$, $\rx/\tip^2$ is bounded by the same constant.\\

Next, let $\tix_+=\tix_+(\rx,\tip)$ and $\tix_- = \tix_-(\rx,\tip)$ be the critical points of $\tix\mapsto h(\tip,\rx,\tix)$ defined in \eqref{def:xstar}. Then $H_0(\rx,\tip,\tik) = \hat{H}(\rx,\tip,\tik)$ is equivalent to
saying that $h(\sqrt{\rx},\rx,\tip,\tik) = h(\tix_+,\rx,\tip,\tik)$, so by the mean-value theorem, there is at least one point $\hat{x} \in (\sqrt{\rx},\tix_+)$ at which $\dfdx{h}{\tix}(\hat{x},\rx,\tip,\tik) = 0$. But the only positive roots of $\tix\mapsto \dfdx{h}{\tix}(\tix,\rx,\tip,\tik)$ are $\tix_-$ and $\tix_+$; therefore, $\hat{x} = \tix_- \in (\sqrt{\rx},\tix_+)$, from which it follows that
$$ \tix_- =  \left(\frac{\rx + \tip - \sqrt{A\rx-B}}{2}\right)^{1/2} > \sqrt{\rx} \quad \implies \quad \rx + \sqrt{A\rx - B} < \tip. $$
Thus, the definition of $\tix_+$ leads to  
$ \tip/2 \leq \tix_+^2 < \tip. $\\

We now rewrite the equation $\gamma(\rx,\tip) = 0$ in \eqref{def:gamma} as 
\begin{equation}
\label{eq:log_equal}
 - 2\sqrt{\rx} + \log\left(\frac{(\tip-\sqrt{\rx})^2}{(\tip+\sqrt{\rx})^2}\right)
=  -2\tix_+ + \log\left(\frac{(\tip-\tix_+)^2+\tix_+^2-\rx}{(\tip+\tix_+)^2+\tix_+^2-\rx}\right).
\end{equation}
Notice that both logarithms are negative or zero, since their arguments are less than or equal to 1. This
allows us to write for instance
$$ 2\tix_+ - 2\sqrt{\rx} + \log\left(\frac{(\tip-\sqrt{\rx})^2}{(\tip+\sqrt{\rx})^2}\right) < 0, $$
or, equivalently
$$ \sqrt{\rx} - \log\left(\frac{1-\frac{\sqrt{\rx}}{\tip}}{1+\frac{\sqrt{\rx}}{\tip}}\right) >  \tix_+.$$
We now recall that the Taylor expansion
$$ -\log\left(\frac{1-a}{1+a}\right) = 2\left(a +\frac{a^3}{3} + \frac{a^5}{5} + \cdots\right)$$
implies
\begin{equation}
\label{eq:log_ineq}
 2a \leq -\log\left(\frac{1-a}{1+a}\right)  \leq \frac{2a}{1-a^2} \qquad \text{for $0\leq a < 1$.} 
\end{equation}
By letting $a = \sqrt{\rx}/\tip < 1/5^{1/4} < 1$, we deduce that
$$ \tix_+ < \sqrt{\rx}\left(1 + \frac{2/\tip}{1-1/\sqrt{5}}\right) = \frac{\sqrt{\rx}}{\tip}\left({\tip} + \frac{5+\sqrt{5}}{2}\right)\leq \frac{C\sqrt{\rx}}{\tip}, $$
where $C = \frac{5+3\sqrt{5}}{2} \approx 5.8541$. Squaring both sides and using the fact that $\tix_+^2 \geq \tip/2$ leads to
$$ \rx \geq \frac{\tip}{2} \cdot \frac{\tip^2}{C^2} = C_1\tip^3,\qquad \text{with $C_1 = \frac{1}{2C^2} = \frac{7-3\sqrt{5}}{20} \approx 0.0146.$} $$
Now, to prove the upper bound $\rx \leq C_2\tip^3$, we rewrite \eqref{eq:log_equal}
as
\begin{equation}
\label{eq:myeq}
\sqrt{\rx} -  \log\left(\frac{1-\frac{\sqrt{\rx}}{\tip}}{1+\frac{\sqrt{\rx}}{\tip}}\right) =\tix_+ -\frac12 \log\left(\frac{(\tip-\tix_+)^2+\tix_+^2-\rx}{(\tip+\tix_+)^2+\tix_+^2-\rx}\right).
\end{equation}
Using \eqref{eq:log_ineq}, we see that the left-hand side of \eqref{eq:myeq} is bounded below by 
\begin{equation}
\label{eq:rx_lowerbound}
\sqrt{\rx} -  \log\left(\frac{1-\frac{\sqrt{\rx}}{\tip}}{1+\frac{\sqrt{\rx}}{\tip}}\right)
\geq \sqrt{\rx} + \frac{2\sqrt{\rx}}{\tip} = \sqrt{\rx}\left(1 + \frac{2}{\tip}\right).
\end{equation}
For the right-hand side of \eqref{eq:myeq}, we observe that
$$\tix_+ -\frac12 \log\left(\frac{(\tip-\tix_+)^2+\tix_+^2-\rx}{(\tip+\tix_+)^2+\tix_+^2-\rx}\right) = 
\tix_+ - \frac12\log\left(\frac{1-\frac{2\tip\tix_+}{\tip^2+2\tix_+^2-\rx}}{1+\frac{2\tip\tix_+}{\tip^2+2\tix_+^2-\rx}}\right). $$
Since $\rx \leq \tip^2/\sqrt{5}$, we have
$$ \tip^2 + 2\tix_+^2 - \rx \geq 0.5528\tip^2 + 2\tix_+^2 \geq 2(0.5528\cdot 2)^{1/2}\tip\tix_+ = 2.1029\tip\tix_+,$$
where we used the algebraic inequality $D_1^2\tip^2 + D_2^2\tix_+^2 \geq 2D_1D_2\tip\tix_+$.
Therefore, we can invoke our Taylor expansion with $a = \frac{2\tip\tix_+}{\tip^2+2\tix_+^2-\rx} < 0.9511 < 1$ to deduce that
\begin{equation}
\label{eq:dx_upperbound}
\tix_+ -\frac12 \log\left(\frac{(\tip-\tix_+)^2+\tix_+^2-\rx}{(\tip+\tix_+)^2+\tix_+^2-\rx}\right) < \tix_+ + \frac{a}{1-(0.9511)^2}\leq \tix_+ + \frac{20.9443 \tip\tix_+}{0.5528\tip^2+2\tix_+^2}.
\end{equation}
Inserting the bounds \eqref{eq:rx_lowerbound} and \eqref{eq:dx_upperbound} into \eqref{eq:myeq}, we deduce that
$$
\sqrt{\rx}\left(1 + \frac{2}{\tip}\right) < 
\tix_+ + \frac{20.9443 \tip\tix_+}{0.5528\tip^2+2\tix_+^2}.
$$
Rearranging leads to
\begin{align*}
 \sqrt{\rx} &\leq \tix_+\left(\frac{{\tip}}{2+{\tip}} + \frac{20.9443\tip^2}{(2+{\tip})(0.5528\tip^2+2\tix_+^2)}\right)\\
& \leq \tix_+\left(\frac{{\tip}}{2+{\tip}} + \frac{20.9443{\tip}}{(2+{\tip})(1+0.5528{\tip})}\right)&& \text{since  $\tix_+^2 \geq \tip/2$,}\\
&\leq 10.9721 {\tip}\tix_+ && \text{since ${\tip} \geq 0$,}\\
&\leq 10.9721 \tip^{3/2} && \text{since $\tix_+ \leq \sqrt{\tip}$.}
\end{align*}
Squaring both sides leads finally to $\rx \leq C_2\tip^3$, with $C_2 = (10.9721)^2 \approx 120.4$.
\bibliographystyle{plain}
\bibliography{References}
\end{document}